

\magnification=\magstep1


\def\res{{\rm Res}} 
\def\P{{\bf P}}     
\def\C{{\bf C}}     
\def\R{{\bf R}}     
\def\Q{{\bf Q}}     
\def\Z{{\bf Z}}     
\def\ox{{\cal O}_X} 
\def\UU{{\cal U}}   
\def\VV{{\cal V}}   
\def\WW{{\cal W}}   
\def\CC{{\cal C}}   
\def\DD{{\cal D}}   
\def\EE{{\cal E}}   
\def\OO{{\cal O}}   

\def\curr#1#2#3{{'\!\DD^{#1,#2}_{#3}}}
\def\cinf{C^{\infty}}
\let\\=\cr\tabskip\centering

\overfullrule=0pt
\def\noi{\noindent}

\def\hfl#1#2{\smash{\mathop{\hbox to 12mm{\rightarrowfill}}
\limits^{\scriptstyle#1}_{\scriptstyle#2}}}

\font\bigbf=cmbx10 scaled \magstep1  

\def\linebreak{\relax\ifhmode\unskip\break\else
    \nonhmodeerr@\linebreak\fi}
\def\newline{\relax\ifhmode\null\hfil\break
\else\nonhmodeerr@\newline\fi}
\let\dsize\displaystyle
\def\frac#1#2{{#1\over#2}}

\def\supp{{\rm supp}}
\def\nw{\eta({\omega})}
\def\ip{\,\smash{\raise7pt\hbox
{\leavevmode\hbox{\vtop{\kern7pt\hrule
  width4pt}\vrule}}}\,\,}  
\def\mapname#1{\ \smash{\mathop{\longrightarrow}\limits^{#1}}\ } 
\def\mapdown#1{\Big\downarrow\rlap{$\vcenter{\hbox{$\scriptstyle#1$}}$}}
\def\dbar{\bar{\partial}} 
\def\efes{F_0,\ldots,F_n}


%
%
%
%

\def\mapdown#1{\Big\downarrow
	 \rlap{$\vcenter{\hbox{$\scriptstyle#1$}}$}}

\def\noe#1{
      \vcenter{\hbox{$\!\scriptstyle#1$}} \nearrow }
\def\soe#1{
      \vcenter{\hbox{$\!\scriptstyle#1$}} \searrow }

\centerline {\bigbf Residues in Toric Varieties}
\vskip .5cm
\centerline {June 22, 1995}
\vskip .7cm
\centerline{{\bf  Eduardo Cattani}}
\centerline{Department of Mathematics and Statistics}
\centerline{University of Massachusetts}
\centerline{Amherst, MA 01003, \ U.S.A. }

\vskip .4cm

\centerline{{\bf  David Cox}}
\centerline{Department of Mathematics and Computer Science}
\centerline{Amherst College}
\centerline{Amherst, MA 01002, \ U.S.A. }

\vskip .4cm

\centerline{\bf Alicia Dickenstein}
\centerline{Departamento de Matem\'atica, F.C.E. y N. }
\centerline{Universidad de Buenos Aires}
\centerline{Ciudad Universitaria - Pabell\'on I }
\centerline{1428 Buenos Aires, \ Argentina }
\vskip 1.5cm

\noi {\bf Introduction}
\bigskip

Toric residues provide a tool for the study of certain homogeneous
ideals of the homogeneous coordinate ring of a toric variety---such as
those appearing in the description of the Hodge structure of their
hypersurfaces [BC].  They were introduced in [C2], where some of their
properties were described in the special case when all of the divisors
involved were linearly equivalent.  The main results of this paper
are: an extension of the Isomorphism Theorem of [C2] to the case of
non-equivalent ample divisors, a global transformation law for toric
residues, and a theorem expressing the toric residue as a sum of local
(Grothendieck) residues.

Let us first establish the notation we will use.  We will assume that
$X$ is a complete toric variety of dimension $n$.  As such, $X$ is
determined by a fan $\Sigma$ in an $n$-dimensional real vector space
$N_{\R}$.  There is a distinguished lattice of maximal rank $N \subset
N_{\R}$ and we let $M$ denote the dual lattice.  The $N$-generators of
the 1-dimensional cones in $\Sigma$ will be denoted
$\eta_1,\ldots,\eta_{n+r}$. This means that $r$ is the rank of the
Chow group $A_{n-1}(X)$.  We will make frequent use of the homogeneous
coordinate ring $S$ of $X$, which is the polynomial ring $S =
\C[x_1,\ldots,x_{n+r}]$.  Here, each variable $x_i$ corresponds to the
generator $\eta_i$ and hence to a torus-invariant irreducible divisor
$D_i$ of $X$.  As in [C1], we grade $S$ by declaring that the monomial
$\Pi_{i=1}^{n+r} x_i^{a_i}$ has degree $[\sum_{i=1}^{n+r} a_i\,D_i]
\in A_{n-1}(X)$.

We will let $\beta = \sum_{i=1}^{n+r} \deg(x_i) \in A_{n-1}(X)$ denote
the sum of the degrees of the variables.  As is well known, $\beta$ is
the anticanonical class on $X$.  Then, given homogeneous polynomials
$F_i \in S_{\alpha_i}$ for $i = 0,\ldots,n$, we define their {\it
critical degree\/} to be
$$
\rho = \big(\textstyle{\sum_{i=0}^n \alpha_i}\big) - \beta \in
A_{n-1}(X).
$$
As in [C2], each $ H \in S_\rho$ determines a meromorphic $n$-form on
$X$
$$
\omega_F(H) = {H\,\Omega \over F_0\cdots F_n}\,,
$$
where $F$ stands for the vector $(F_0,\ldots,F_n)$
and $\Omega$ is a choice of an Euler form in $X$ [BC].
If the $F_i$ don't
vanish simultaneously on $X$, then relative to the open cover $U_i =
\{x \in X : F_i(x) \ne 0\}$ of $X$, this gives a \v Cech cohomology
class $[\omega_F(H)] \in H^n(X,\widehat\Omega_X^n)$.  Here,
$\widehat\Omega_X^n$ is the sheaf of Zariski $n$-forms on $X$, i.e.,
$\widehat\Omega^n_X = j_*\Omega^n_{X_0}$, where $X_0$ is the smooth
part of $X$ and $j\colon X_0\hookrightarrow X$ is the natural
inclusion.

It is not hard to see that $[\omega_F(H)]$ is alternating on the order
of $\efes$ and that if $H$ is in the ideal $\langle
F_0,\ldots,F_n\rangle$, then $\omega_F(H)$ is a \v Cech coboundary.
Thus, $[\omega_F(H)]$ depends only on the equivalence class of $H$
modulo the ideal generated by $\efes$.  Then the {\it toric residue\/}
$$
\res_F : S_\rho/\langle F_0,\ldots,F_n\rangle_\rho \longrightarrow \C
$$
is given by the formula
$$
\res_F(H) = {\rm Tr}_X([\omega_F(H)]),
$$
where ${\rm Tr}_X : H^n(X,\widehat\Omega_X^n) \to \C$ is the trace map.
When there is no danger of confusion, we will write $\res(H)$ instead of
$\res_F(H)$.
\smallskip

Our first main result is the following {\it Global Transformation
Law}.

\proclaim Theorem 0.1. Let $F_i \in S_{\alpha_i}$ and $G_i \in
S_{\beta_i}$ for $i = 0,\ldots,n$.  Suppose
$$
G_j = \sum_{i=0}^n A_{ij}\,F_i,
$$
where $A_{ij}$ is homogeneous of degree $\beta_j - \alpha_i$, and
assume the $G_i$ don't vanish simultaneously on $X$.  Let $\rho$ be
the critical degree for $F_0,\ldots,F_n$.  Then, for each $H \in
S_\rho$, $H\det(A_{ij})$ is of the critical degree for
$G_0,\ldots,G_n$, and
$$
\res_F(H) = \res_G(H\det(A_{ij})).
$$

The proof uses a \v Cech cochain argument.  One application of this
transformation law is that in certain cases, we can describe explicit
elements of $S_\rho$ with nonzero residue.  For this purpose, assume
$X$ is complete and its fan $\Sigma$ contains a $n$-dimensional {\it
simplicial} cone $\sigma$.  Then denote the variables of the
coordinate ring as $x_1,\ldots,x_n,z_1,\ldots,z_r$, where
$x_1,\ldots,x_n$ correspond to the 1-dimensional cones of $\sigma$.
Also suppose that $\alpha_0,\ldots,\alpha_n$ are ${\bf Q}$-ample
classes, which means that some multiple is Cartier and ample.  In this
situation, we will show that each $F_j \in S_{\alpha_j}$ can be
written in the form
$$
F_j = A_{0j}\,z_1\cdots z_r + \sum_{i=1}^n A_{ij}\, x_i.
$$
Then the $(n+1)\times(n+1)$-determinant $\Delta_\sigma = \det(A_{ij})$
is in $S_\rho$ and has the following important property.

\proclaim Theorem 0.2. Assume $X$ is complete and $\sigma \in \Sigma$
is simplicial and $n$-dimensional.  Suppose that $F_i\in S_{\alpha_i}$
for $i = 0,\ldots,n$, where $\alpha_i$ is $\Q$-ample and the $F_i$
don't vanish simultaneously on $X$.  Then
$$
\res_F(\Delta_\sigma) = \pm1.
$$

The Global Transformation Law allows us to reduce the proof of this
theorem to the special case when $F_0 = z_1\cdots z_r$ and $F_i=x_i$,
$i=1,\ldots,n$.  This is done in \S 2.  An alternate proof for
simplicial toric varieties is given in \S 4 as an application of
Theorem~0.4.
\medskip

In \S3 we prove the following {\it Residue Isomorphism Theorem\/}.

\proclaim Theorem 0.3. Let $X$ be complete and simplicial,
and assume that $F_i\in S_{\alpha_i}$ for $i = 0,\ldots,n$, where
$\alpha_i$ is ample and the $F_i$ don't vanish simultaneously on $X$.
Then:
\vskip0pt
\item{(i)} The toric residue map $\res_F : S_\rho/\langle
F_0,\ldots,F_n\rangle_\rho \to \C$ is an isomorphism.
\item{(ii)} For each variable $x_i$, $0 \le i \le n+r$, we have
$x_i\cdot S_\rho \subset \langle F_0,\ldots,F_n\rangle$.
\vskip0pt

In the case when all the $\alpha_i$ are equal to a fixed ample divisor
$\alpha$, this theorem follows from the fact that $F_0,\ldots, F_n$
are a regular sequence in the Cohen-Macaulay ring $S_{*\alpha} =
\oplus_{k\geq 0}S_{k\alpha}$ [C2, \S3]. In the general case, the proof
relies on the use of the Cayley trick and results of Batyrev and Cox
[BC] concerning the cohomology of projective hypersurfaces in toric
varieties, to show that
$$
\dim_\C(S_\rho/\langle F_0,\ldots,F_n\rangle_\rho) = 1
$$
when $X$ is simplicial and the divisors $F_i = 0$ are ample with empty
intersection. Then, the first (and main) part of the Residue
Isomorphism Theorem follows immediately from Theorem~0.2, and the
second part is a consequence of the first using Theorem 0.2 and
Cramer's Rule.

As a corollary of Theorems~0.2 and~0.3, we get a simple algorithm (see
Remark 3.11) for computing toric residues in terms of normal forms.
In \S3, we will also conjecture a more general form of Theorem 0.3 and
give some examples in support.
\medskip

The main result of \S4 is a theorem stating that for simplicial toric
varieties, the toric residue may be computed as a sum of local
Grothendieck residues.  The toric setting is not essential here and,
in fact, it is convenient to work with the more general notion of a
$V$-manifold or orbifold.  The proof of the following local/global
theorem is based on the theory of residual currents ([CH]).

\proclaim Theorem 0.4. Let $X$ be a complete simplicial toric
variety of dimension $n$, and let $F_0,\ldots,F_n$ be homogeneous
polynomials which don't vanish simultaneously on $X$.  If $H$ is a polynomial
in $S_\rho$, where $\rho$ is the critical degree, and $D_{\hat k} = \{x \in
X: F_i(x) = 0,\ i \ne k\}$ is finite, then the toric residue is given
by
$$
\res_F(H) =
(-1)^{k}\, \sum_{x\in D_{\hat k} }\res_{k,x}\biggl({H\,\Omega \over
F_0\cdots F_n}\biggr).
$$

Here, as we will explain in \S4, $\dsize{\res_{k,x}\Bigl({H\,\Omega
\over F_0\cdots F_n}\Bigr)}$ denotes the local Grothendieck residue
$\dsize{\res_{x}\Bigl({(H/F_k)\,\Omega \over
F_0\cdots\widehat{F_k}\cdots F_n}\Bigr)}$.  Note that the finiteness
condition holds automatically whenever the divisor $\{F_k = 0\}$ is
$\Q$-ample.  Under appropriate conditions, Theorem 0.4 gives a
framework for the study of sums of local residues---both in the affine
and toric cases---as a global residue defined in a suitable toric
compactification.  It is possible, for example, to interpret in this
light the results of [CDS] which correspond to the case when the toric
variety under consideration is a weighted projective space.
\medskip

Finally, in \S5, we show that, in the equal degree case, the toric
residue equals a single local residue at the origin of the affine cone
of $X$.  This generalizes the observation in [PS] that toric residues
on $\P^n$ can be written as a residue at the origin in $\C^{n+1}$.
\bigskip

\noi {\bf Acknowledgements.} \ Eduardo Cattani was supported
by NSF Grant DMS-9404642; part of the work on this paper was done while
he was visiting the Institut Fourier, Grenoble, and the University of
Buenos Aires. He is grateful for their support and hospitality. David Cox
was supported by NSF Grant DMS-9301161. Alicia Dickenstein was supported
by UBACYT and CONICET, Argentina.
\bigskip
\bigskip

\noi {\bf \S1. The Global Transformation Law}
\bigskip

This section will prove the Global Transformation Law (Theorem 0.1)
for toric residues on an arbitrary $n$-dimensional complete toric
variety $X$.  Given $F_i$ and $G_j = \sum_{i=0}^n A_{ij}\,F_i$ as in
the statement of the theorem, first observe that if $G_0,\ldots,G_n$
have no common zeroes in $X$, then the same holds for
$F_0,\ldots,F_n$.  Thus we get open covers $U_i = \{x\in X : F_i(x)
\ne 0\}$ and $V_j = \{x \in X : G_j(x) \ne 0\}$ of $X$, which we
denote ${\cal U}$ and ${\cal V}$ respectively.

If the critical degrees of the $F_i$ and $G_j$ are $\rho(F)$ and
$\rho(G)$ respectively, then
$$
\rho(G) = \rho(F) + \deg\bigl(\det(A_{ij})\bigr)
$$
follows easily since $A_{ij}$ homogeneous and $\deg A_{ij} = \deg G_j
- \deg F_i$.  Thus, if $H$ is homogeneous of degree $\rho(F)$, then
we get \v Cech cohomology classes $[\omega_F(H)] \in H^n({\cal
U},\widehat\Omega_X^n)$ and $[\omega_G(H\det(A_{ij}))] \in H^n({\cal V},
\widehat\Omega_X^n)$.  To prove Theorem 0.1, it suffices to show that
these cohomology classes have the same image in
$H^n(X,\widehat\Omega_X^n)$.

Consider the open covering ${\cal W} = {\cal U}\cup{\cal V}$.  Since
${\cal U}$ and ${\cal V}$ can be regarded as refinements of ${\cal W}$
with obvious refinement maps, we get a commutative diagram
$$
\matrix{&& H^n({\cal U},\widehat\Omega_X^n) &&\cr
& \nearrow && \searrow & \cr
H^n({\cal W},\widehat\Omega_X^n) &&&& H^n(X,\widehat\Omega_X^n).\cr
& \searrow && \nearrow &\cr
&& H^n({\cal V},\widehat\Omega_X^n) &&\cr}
$$
Then Theorem 0.1 is an immediate consequence of the following
proposition.

\proclaim Proposition 1.1. There is a cohomology class $[\theta] \in
H^n({\cal W},\widehat\Omega_X^n)$ which maps to both $[\omega_F(H)]
\in H^n({\cal U},\widehat\Omega_X^n)$ and $[\omega_G(H\det(A_{ij}))]
\in H^n({\cal V},\widehat\Omega_X^n)$ in the above diagram.

\noi{\bf Proof.}\ \ We first introduce some notation for the \v Cech
complex of $\WW = \UU \cup \VV$.  Given index sets $I=\{0\leq
i_1<\cdots<i_p \leq n\}$ and $J=\{0\leq j_1<\cdots<j_q \leq n\}$ with
$p = |I|$ and $q = |J|$, let $W_{IJ} = \bigcap_{i\in I} U_i\cap
\bigcap_{j\in J} V_j$.  Also, let $I'$ denote the complement of $I$ in
$\{0,\ldots,n\}$, ordered increasingly.

We define a \v Cech cochain $\theta \in \CC^n(\WW,\widehat\Omega^n_X)$
by the assignment
$$
W_{IJ} \mapsto \theta_{IJ} = \epsilon(I)\ {{H\,\det(M_{I'J})\,
\Omega}\over{F_I\,G_J}}\in \widehat\Omega^n_X(W_{IJ}).
$$
Here, $I$ and $J$ are index sets with $p+q = n+1$, $M_{I'J}$
is the $I'J$-minor of the matrix $(A_{ij})$, $\epsilon(I)$ is the
sign of the permutation $(I,I')$, $F_I=F_{i_1}\cdots F_{i_p}$, and
$G_J=G_{j_1}\cdots G_{j_q}$.

When $I = \{0,\ldots,n\}$, both $I'$ and $J$ are empty and
$\theta$ assigns to the open set $U_0\cap\cdots\cap U_n$ the form
$H\,\Omega/(F_0\cdots F_n)$.  Thus, the refinement map
$\CC^n(\WW,\widehat\Omega^n_X) \to \CC^n(\UU,\widehat\Omega^n_X)$ maps
$\theta$ to the cocycle $\omega_F(H)$.  Similarly, when $J =
\{0,\ldots,n\}$, the cochain $\theta$ assigns to the open set
$V_0\cap\cdots\cap V_n$ the form $H\,\det(A_{ij})\,\Omega/(G_0\cdots
G_n)$ and, hence, the refinement map $\CC^n(\WW,\widehat\Omega^n_X)
\to \CC^n(\VV,\widehat\Omega^n_X)$ maps $\theta$ to the cocycle
$\omega_G(H\,\det(A_{ij}))$.  Consequently, the proposition will
follow once we show that $\theta$ is also a cocycle, i.e.,
$\delta(\theta) = 0$, where $\delta : \CC^n({\cal
W},\widehat\Omega_X^n) \to \CC^{n+1}({\cal W},\widehat\Omega_X^n)$ is
the \v Cech coboundary.

To prove that $\delta(\theta) = 0$, let $I$ and $J$ be index sets with
$p+q=n+2$.  Then set $I_k = I - \{i_k\}$, $1\leq k \leq p$,
$I_k'=I'\cup \{i_k\}$, and $J_{\ell} = J - \{j_{\ell}\}$, $1\leq
{\ell} \leq q$, so that
$$
\eqalign{
(\delta\theta)_{IJ} = &\ \sum_{k=1}^p\ (-1)^{k-1}\ \theta_{I_kJ}
+ (-1)^p \sum_{\ell=1}^q\ (-1)^{\ell-1}\ \theta_{IJ_{\ell}}\cr
= &\ \sum_{k=1}^p\ \frac{(-1)^{k-1}\ \epsilon(I_k) \
\det(M_{I_k'J})\ H\ \Omega}{F_{I_k}\ G_J}\cr
&\ + \sum_{\ell=1}^q\ \frac{(-1)^{p + \ell-1}\
\epsilon(I)\ \det(M_{I'J_{\ell}})\ H\ \Omega}{F_{I}\
G_{J_{\ell}}}.\cr}
$$
Writing the last expression with common denominator
$F_I\,G_{J}$, it suffices to show that
$$
\sum_{k=1}^p\ {(-1)^{k-1}\, \epsilon(I_k) \ F_{i_k}\,
\det(M_{I_k'J}})
+ \sum_{\ell=1}^q\ {(-1)^{p + \ell-1}\,
\epsilon(I)\ G_{j_{\ell}}\, \det(M_{I'J_{\ell}})} = 0.
$$
If we substitute $G_{j_{\ell}}$ by $ \sum_{k=1}^p\, A_{i_kj_{\ell}}\,
F_{i_k} + \sum_{u\in I'} \, A_{uj_{\ell}}\, F_{u}$, then the above
equation becomes
$$
\eqalign{\sum_{k=1}^p\, \left[{(-1)^{k-1}\, \epsilon(I_k) \,
\det(M_{I_k'J})}
+ \sum_{\ell=1}^q\, {(-1)^{p + \ell-1}\,
\epsilon(I)\ A_{i_kj_{\ell}}\, \det(M_{I'J_{\ell}})}\right]\,
F_{i_k}&\cr
\null + \epsilon(I)\, (-1)^p\, \sum_{u\in I'}\left[
\sum_{\ell=1}^q\, {(-1)^{ \ell-1}\,
A_{uj_{\ell}}\, \det(M_{I'J_{\ell}})}\right]\, F_u& = 0.
\cr}
$$
We will show that the expressions inside the two sets of brackets are
identically zero.

First, for $u\in I'$, note that $\sum_{\ell=1}^q\, {(-1)^{\ell-1}\,
A_{uj_{\ell}}\det(M_{I'J_{\ell}})}$ is the determinant of the matrix
whose first row consists of $A_{uj_{\ell}}$, $\ell = 1,\ldots,q$ and
whose remaining rows are the same as those of the matrix $M_{I'J}$.
But, since $u\in I'$, such a matrix has two identical rows and its
determinant vanishes.  Hence the sum inside the second set of brackets
is zero.

Turning to the first set of brackets, note that expanding the
determinant of the $q\times q$-matrix $M_{I_k'J}$ along the row
corresponding to $i_k$ gives
$$
\det(M_{I_k'J}) = (-1)^m\, \sum_{\ell=1}^q\, {(-1)^{ \ell-1}\,
A_{i_kj_{\ell}}\, \det(M_{I'J_{\ell}})},
$$
where $m$ denotes the number of elements of $I'$ which precede $i_k$.
However, since going from $(I,I')$ to $(I_k,I_k')$ requires $(m-k+p)$
transpositions, we see that
$$\epsilon(I_k) = (-1)^{m-k+p}\ \epsilon(I),$$
and it follows that the desired expression is zero.$\ \ \diamond$
\bigskip
\bigskip

\noi {\bf \S2. Elements with Nonzero Residue}
\bigskip

The goal of this section is to prove Theorem 0.2.  We begin with $F_j
\in S_{\alpha_j}$, $0 \le j \le n$, which don't vanish simultaneously
on the complete toric variety $X$.  As in the introduction, we have
the coordinate ring $S = {\bf C}[x_1,\ldots,x_{n+r}]$ of $X$, where
the variables $x_i$ correspond to the ${\bf Z}$-generators $\eta_i$ of
the 1-dimensional cones of $\Sigma$.  For each $n$-dimensional cone
$\sigma\in \Sigma$, set $\hat x_{\sigma} =
\prod_{\eta_i\not\in\sigma}\,x_i$ and let $B(\Sigma)
\subset S$ be the ideal generated by the monomials $\hat x_{\sigma}$,
$\sigma\in\Sigma$.

We also assume that $\alpha_j$ is ${\bf Q}$-ample.  This means that
$d\,\alpha_j$ is ample for some positive integer $d$, so that
$S_{d\alpha_j} \subset B(\Sigma)$ by Lemma 9.15 of [BC].  Then
$(S_{\alpha_j})^d \subset S_{d\alpha_j} \subset B(\Sigma)$, and since
$B(\Sigma)$ is radical, we conclude that
$$
S_{\alpha_j} \subset B(\Sigma)\leqno{(2.1)}
$$
when $\alpha_j$ is ${\bf Q}$-ample.

To see the relevance of (2.1), fix a $n$-dimensional cone $\sigma \in
\Sigma$ and, as in Theorem 0.2, assume that $\sigma$ is {\it
simplicial\/}.  Then we can arrange for $\eta_1,\ldots,\eta_n$ to be
the generators of $\sigma$, and we make a slight notational change
replacing the variable $x_{n+a}$ by $z_a$, $a=1,\ldots,r$.  Then each
$F_j$ can be written
$$
F_j = B_j + \sum_{i=1}^n A_{ij}\,x_i,
$$
where $B_j$ depends only on $z_1,\ldots,z_r$.  But $F_j \in B(\Sigma)$
by (2.1) and, since $B(\Sigma)$ is a monomial ideal, it follows easily
that $B_j$ must be divisible by $\hat x_{\sigma} = z_1\cdots z_r$.
Thus $B_j = A_{0j}\,z_1\cdots z_r$, so that
$$
F_j = A_{0j}\,z_1\cdots z_r + \sum_{i=1}^n A_{ij}\,x_i,
\leqno{(2.2)}
$$
as claimed in the introduction.  Then we define
$$
\Delta_{\sigma} = \det(A_{ij}),\quad 0\leq i,j \leq n.\leqno{(2.3)}
$$
Note that
$$\deg(\Delta_{\sigma}) = \sum_{j=0}^n \alpha_i - \sum_{i=1}^n
\deg(x_i) - \sum_{a=1}^r
\deg(z_a) = \rho(F)\,.$$

A direct application of the Global Transformation Law
to (2.2) and (2.3) yields
$$
\res\biggl(\frac{\Delta_{\sigma}\,\Omega}{F_0\cdots F_n} \biggr) =
\res\biggl(\frac{\Omega}{(z_1\cdots z_r)\cdot x_1\cdots x_n}\biggr).
$$
Thus, to prove Theorem 0.2, we need only prove the following proposition.

\proclaim Proposition 2.4. With assumptions and notations as
above,
$$
\res\biggl(\frac{\Omega}{(z_1\cdots z_r)\cdot x_1\cdots x_n}\biggr)
= \pm 1.
$$

\noi{\bf Proof.}\ \  We first study the restriction of
$$
\omega = {\Omega\over(z_1\cdots z_r)\cdot x_1\cdots x_n}\leqno{(2.5)}
$$
to the affine open set $X_\sigma \subset X$ determined by $\sigma$.
To apply the construction of [C1] to the simplicial toric variety
$X_\sigma$, we start with the exact sequence
$$
0 \longrightarrow M \mapname{\gamma} \Z^n \longrightarrow D(\sigma)
\longrightarrow 0,
$$
where $\gamma(m) = (\langle m,\eta_1\rangle,\ldots,\langle
m,\eta_n\rangle)$.  Then $D(\sigma)$ is finite since $\sigma$ is
simplicial, and $G(\sigma) = {\rm Hom}_\Z(D(\sigma),\C^*)$ is
isomorphic to $N/N'$, where $N'$ is the sublattice of $N$ generated by
$\eta_1,\ldots,\eta_n$.  The map $\Z^n \to D(\sigma)$ induces an
action of $G(\sigma)$ on $\C^n$, and it follows from [C1] that the
quotient is $X_\sigma$.  In particular, we have a quotient map
$p_\sigma : \C^n \to X_\sigma$.

To relate this construction to $S = \C[x_1,\ldots,x_{n+r}] =
\C[x_1,\ldots,x_n,z_1,\ldots,z_r]$, note that $S$ is the coordinate
ring of the affine space $\C^{n+r}$.  Then let $Z(\Sigma) \subset {\bf
C}^{n+r}$ be the subvariety defined by $B(\Sigma)$.  In [C1], a
natural map $p:{\bf C}^{n+r} - Z(\Sigma) \to X$ is constructed.  Now
consider the inclusion ${\bf C}^n \to {\bf C}^{n+r}$ defined by
$$
(x_1,\ldots,x_n) \mapsto (x_1,\ldots,x_n,1,\ldots,1). \leqno{(2.6)}
$$
One easily sees that the image of this map lies in $\C^{n+r} -
Z(\Sigma)$, and the argument of Theorem 1.9 of [BC] shows that we have
a commutative diagram
$$
\matrix{\C^n & \mapname{p_\sigma} & X_\sigma\cr
	\downarrow && \downarrow \cr
	\C^{n+r} - Z(\Sigma) & \mapname{p} & X.\cr} \leqno{(2.7)}
$$

We now recall the Euler form $\Omega$ from [BC].  Fix an integer basis
$m_1,\ldots,m_n$ for the lattice $M$.  Then, given a subset $I =
\{\ell_1,\ldots,\ell_n\} \subset \{1,\ldots,n+r\}$ with $n$ elements,
define
$$ \det(\eta_I) = \det(\langle m_i,\eta_{\ell_j}
	\rangle_{\scriptscriptstyle{1 \le i,j \le n}}).
$$
Also set $dx_I = dx_{\ell_1}\wedge\cdots\wedge dx_{\ell_n}$ and
$\hat{x}_I = \Pi_{\ell \notin I} x_\ell$.  Then $\Omega$ is given by
the formula
$$
\Omega = \sum_{|I| = n} \det(\eta_I)\, \hat{x}_I\, dx_I,\leqno{(2.8)}
$$
where the sum is over all $n$-element subsets $I \subset
\{1,\ldots,n+r\}$.  Note that $\Omega$ is well-defined up to $\pm1$.

{}From (2.6) and (2.7), we see that $p_\sigma^*(\omega)$ is computed by
setting $z_1 = \cdots = z_r = 1$ in the above formula for $\Omega$.
Thus
$$
p_\sigma^*(\omega) = \pm{|N/N'|\, dx_1\wedge\cdots\wedge dx_n \over
x_1\cdots x_n}. \leqno{(2.9)}
$$

The next step in the proof is to study what happens when we change $X$
slightly.  Suppose that $\Sigma'$ is a refinement of the fan $\Sigma$
such that $\sigma$ is still a cone of $\Sigma'$.  Then we get a
birational morphism $\pi:X' \to X$ of toric varieties which is the
identity on the affine piece $X'_\sigma = X_\sigma$.  For $X'$, we
have an Euler form $\Omega'$, and the analog of $\omega$ in (2.5) is
denoted $\omega'$.  Note that $x_1,\ldots,x_n$ have the same meaning
for both $\omega$ and $\omega'$.  Then, since $\sigma$ is a cone for
both fans and we haven't changed $N$, it follows from (2.9) that
$\pi^*(\omega) = \omega'$.  We also have open covers $\UU$ of $X$ and
$\UU'$ of $X'$, and we leave it to the reader to verify that these
covers are compatible under $\pi$.  It follows that
$$
\pi^*([\omega]) = [\omega']
$$
as elements
of $H^n(X',\widehat\Omega_{X'}^n)$.  Since
$$
{\rm Tr}_{X'}\circ \pi^* = {\rm Tr}_X
$$
when $\pi$ is birational, $\omega$ and $\omega'$ have the same toric
residue.  In particular,
$$
\res(\omega) = \pm1 \quad \iff \quad \res(\omega') = \pm1.\leqno{(2.10)}
$$

Next, instead of changing the fan, suppose we change the lattice.  If
$N_1$ is a sublattice of $N$, then as explained in [O, Corollary
1.16], we get a toric variety $X_1$ such that $N/N_1$ acts on $X_1$
with $X$ as quotient.  Let $\pi_1 : X_1 \to X$ be the quotient map.
The toric varieties $X$ and $X_1$ have the same coordinate ring $S$
(though the gradings differ).  Now assume that
$\eta_1,\ldots,\eta_{n+r}$ lie in $N_1$.  Then one easily sees that
the Euler forms $\Omega$ and $\Omega_1$ are related by the formula
$$
\Omega = \pm |N/N_1|\,\Omega_1,
$$
so that if $\omega_1$ is the analog of $\omega$ for $X_1$, we have
$$
\pi_1^*(\omega) = \pm |N/N_1|\,\omega_1.
$$
However, since $\pi_1 : X_1 \to X$ is a finite map of degree $|N/N_1|$,
we also have
$$
{\rm Tr}_{X_1} \circ \pi_1^* = |N/N_1|\,{\rm Tr}_X.
$$
{}From here, it follows that $\omega$ and $\omega_1$ have the same toric
residue (up to $\pm1$), and hence
$$
\res(\omega) = \pm1 \quad \iff \quad \res(\omega_1) = \pm1.\leqno{(2.11)}
$$

We can now prove the proposition.  Define $\eta_0 = -\sum_{i=1}^n
\eta_i$, where $\eta_1,\ldots,\eta_n$ are the generators of $\sigma$,
and let $\Sigma'$ be the fan whose cones are generated by proper
subsets of $\{\eta_0,\ldots,\eta_n\}$.  This gives a toric variety
$X'$, and note the $\sigma$ is a cone of $\Sigma'$.  Now let
$\Sigma''$ be the fan consisting of all intersections
$\sigma_1\cap\sigma_2$ for $\sigma_1 \in \Sigma$ and $\sigma_2 \in
\Sigma'$.  Hence $\Sigma''$ is a common refinement of $\Sigma$ and
$\Sigma'$, and $\sigma$ is a cone in all three fans.  The
corresponding toric variety $X''$ maps to both $X$ and $X'$.  Finally,
let $N' \subset N$ be the sublattice generated by
$\eta_1,\ldots,\eta_n$.  Then $\eta_0 \in N'$, and the toric variety
determined by $N'$ and $\Sigma'$ is $\P^n$.  Putting this all
together, we get a diagram of toric varieties
$$
\matrix{&& X'' &&&& \P^n \cr & \swarrow && \searrow && \swarrow \cr
	X &&&& X' &&\cr}
$$
where the first two maps come from refinements which preserve $\sigma$
and the third comes from a change of lattice which preserves the
generators $\eta_i$.  It follows from (2.10) and (2.11) that
$$
\res(\omega) = \pm1
$$
if and only if the corresponding statement is true for $\P^n$.  The
latter is well known; for example, it follows from the Trace Property
for $\P^n$ stated in the introduction to [C2].$\ \ \diamond$
\bigskip

\noi {\bf Remarks 2.12.}\ (i)\ In \S4, we will use the relation between
toric residues and local residues to give a second proof of
Proposition 2.4 in the special case when $X$ is simplicial.
\medskip
\noi (ii)\ If we replace the hypothesis that the $\alpha_i$ are
$\Q$-ample with the weaker assumption (2.1), then the conclusion of
Theorem 0.2 is still true.  In fact, all we need to assume is that
$F_i \in B(\Sigma)$ for all $i$.  This will be useful in \S3.
\medskip
\noi (iii)\ The element $\Delta_{\sigma}$ depends on the choice of
simplicial cone $\sigma$ of dimension $n$ and on the choice of
coefficients $A_{ij}$ in (2.2).  Once Theorem 0.3 is established, it
will follow from Theorem 0.2 that when $X$ is simplicial and the
$\alpha_i$ are ample, the class of $\Delta_{\sigma}$ is unique up to
sign modulo the ideal $\langle F_0,\ldots,F_n\rangle$.  Moreover, if
we pick a basis of $M$ such that $\det(\langle m_i,\eta_j\rangle) >
0$, then one can check that $\res_F(\Delta_{\sigma})=1$.
\medskip
\noi (iv)\  Given {\it any\/} decomposition
$$
F_j = A_{0j}\,z_1\ldots z_r + \textstyle{\sum_{i=1}^n} A_{ij}\,x_i,
$$
the polynomial $\Delta_\sigma = \det(A_{ij})$ satisfies
$\res_F(\Delta_\sigma)=\pm 1$.
\medskip
\noi (v)\ The definition of $\Delta_\sigma$ given by (2.3) generalizes
a well-known construction in the algebraic setting corresponding to
projective space [KK].  Regarding $\P^n$ as a toric variety, we get the
usual graded ring $S=\C[x_0,\ldots,x_n]$, and the ideal $B(\Sigma)$ is
the maximal ideal $\langle x_0,\ldots,x_n\rangle$.  Given
homogeneous polynomials $\efes$, $\deg(F_j) = d_j>0$, whose only
common zero is the origin, let $\sigma$ be the cone whose generators
correspond to the variables $x_1,\ldots,x_n$.  Then
$$
F_j = \frac{1}{d_j} \frac{\partial F_j}{\partial x_0}\,x_0 +
\sum_{i=1}^n \frac{1}{d_j} \frac{\partial F_j}{\partial x_i}
\,x_i
$$
and, consequently, a choice of $\Delta_\sigma$ is given by
$$
\Delta_\sigma = \frac{1}{d_0\cdots d_n} \det\biggl(
\frac{\partial F_i}{\partial x_j}\biggr).
$$
\bigskip

\noindent {\bf \S3. The Codimension One and Residue Isomorphism
Theorems}
\bigskip

Before we can prove the main results of this section, we need to
discuss a toric version of the Cayley trick.  Let $X$ be a complete
toric variety, and let ${\cal L}_0,\ldots,{\cal L}_n$ be ample line
bundles on $X$.  Then consider
$$
Y = \P({\cal L}_0\oplus\cdots\oplus{\cal L}_n).
$$
This variety has a canonical line bundle ${\cal O}_Y(1)$, which is
ample since the ${\cal L}_j$ are ample (see \S1 of Chapter III of
[H]).  Note that $Y$ is a $\P^n$-bundle over $X$, so that $Y$ has
dimension $2n$.

For our purposes, we need to understand $Y$ as a toric variety.  We
begin with the description given in [BB].  The ample line bundle
${\cal L}_j$ is associated to a $n$-dimensional polytope $\Delta_j
\subset M_\R$ of the form
$$
\Delta_j = \{m \in M_\R: \langle m,\eta_i\rangle \ge -a_{ij},\ i =
1,\ldots,n+r\}. \leqno{(3.1)}
$$
Since each ${\cal L}_j$ is ample, the {\it facets\/} of $\Delta_j$
(faces of codimension 1) correspond bijectively to the $\eta_i$, where
$\eta_i$ gives the facet $F_{ij} = \{m \in \Delta_j: \langle
m,\eta_i\rangle = -a_{ij}\} \subset \Delta_j$.

Now consider $\R^n\oplus M_\R$ with the integer lattice $\Z^n\oplus
M$.  Elements of $\R^n\oplus M_\R$ can be uniquely written $\widetilde
m = \lambda_1e_1 + \cdots + \lambda_ne_n + m$, where $\lambda_j \in
\R$ and $m \in M_\R$.  We also have the dual $\R^n\oplus N_\R$ with
lattice $\Z^n\oplus N$, and elements here are written similarly.  Then
define $\Delta \subset \R^n\oplus M_\R$ to be the convex hull
$$
\eqalign{\Delta &= {\rm Conv}\big((\{0\}\times \Delta_0) \cup
(\{e_1\}\times\Delta_1)
\cup \cdots \cup (\{e_n\}\times\Delta_n)\big)\cr
&= \{\lambda_1e_1+\cdots+\lambda_ne_n + \lambda_0 m_0 + \cdots +
\lambda_n m_n : \lambda_j \ge 0,\ \textstyle{\sum_{j=0}^n} \lambda_j =
1,\ m_j \in \Delta_j\}.\cr}\leqno{(3.2)}
$$
This is easily seen to be equivalent to the polytope $\Delta_\sigma$
in Section 3 of [BB].  Since ${\cal O}_Y(1)$ is ample, Proposition 3.2
of [BB] implies that $Y$ is the toric variety determined by the
polytope $\Delta$.  The corresponding fan in $\R^n\oplus N_\R$ is
called the {\it normal fan\/} of $\Delta$.

We next show that the generators of the 1-dimensional cones in the
normal fan are given by
$$
\eqalign{\tilde\eta_i &= \textstyle{\sum_{j=0}^n}
(a_{ij}-a_{i0})e_j + \eta_i,\quad i = 1,\ldots,n+r\cr
	\tilde\nu_0 &= -e_1-\cdots-e_n\cr
	\tilde\nu_j &= e_j,\quad j = 1,\ldots,n.\cr}
$$
The first step is to prove that $\Delta$ is defined by the
inequalities
$$
\eqalign
{\langle \widetilde m,\tilde\eta_i\rangle &\ge -a_{i0}, \quad i =
	1,\ldots,n+r\cr
	\langle \widetilde m,\tilde\nu_0\rangle &\ge -1\cr
	\langle \widetilde m,\tilde\nu_j\rangle &\ge \ 0,\quad j =
	1,\ldots,n.\cr} \leqno{(3.3)}
$$
Write $\widetilde m = \lambda_1e_1+\cdots+\lambda_ne_n + m$, where $m
\in M_\R$, and let $\lambda _0 = 1 - \sum_{j=1}^n\lambda_j$.  Then the
above inequalities are equivalent to
$$
\langle m,\eta_i\rangle \ge -\textstyle{\sum_{j=0}^n}\lambda_j
	a_{ij},\ i = 1,\ldots,n+r,\ \ \
	\lambda_j \ge 0,\ j = 0,\ldots,n,\ \ \
	\textstyle{\sum_{j=0}^n} \lambda_j = 1.
$$

If $\widetilde m \in \Delta$, then (3.2) shows that $m = \sum_{j=0}^n
\lambda_j m_j$, where $m_j \in \Delta_j$, and it follows easily from
(3.1) that the above conditions are satisfied.  Conversely, if
$\widetilde m = \lambda_1e_1+\cdots+\lambda_ne_n + m$ satisfies (3.3),
consider the linear maps $B(m_0,\ldots,m_n) = (-\langle
m_j,\eta_i\rangle)$ and $D(m_0,\ldots,m_n) = \sum_{j=0}^n \lambda_j
m_j$.  Then the above inequalities and the Farkas Lemma (as stated in
Exercise 1.6 of [Z] with $A = C = 0$) imply that there exists
$(m_0,\ldots,m_n)$ with $B(m_0,\ldots,m_n) \le (a_{ij})$ and
$D(m_0,\ldots,m_n) = m$.  This shows that $m = \sum_{j=0}^n \lambda_j
m_j$, where $m_j \in \Delta_j$ by (3.1), and $\widetilde m \in \Delta$
follows immediately.

 From the inequalities defining $\Delta$, we can read off the facets of
$\Delta$ as follows.  First, one easily shows that
$$
\{\widetilde m \in \Delta : \langle \widetilde m,\tilde\eta_i\rangle =
-a_{i0}\}
$$
is the convex hull
$$
{\rm Conv}\big((\{0\}\times F_{i0})\cup(\{e_1\}\times
F_{i1}) \cup\cdots\cup (\{e_n\}\times F_{in})\big),
$$
where $F_{ij} \subset \Delta_j$ is the facet defined by $\eta_i$.
Since this set has dimension $2n-1$, it is a facet of $\Delta$.
Similarly, if one looks at the subsets of $\Delta$ defined by the
equations $\langle \widetilde m,\tilde\nu_0\rangle = -1$ or $\langle
\widetilde m,\tilde\nu_j\rangle = 0$, $1 \le j \le n$, then one gets
the $(2n-1)$-dimensional polytopes obtained by taking the convex hulls
of $n$ of the $n+1$ sets $\{0\}\times\Delta_0,\{e_1\}\times\Delta_1,
\ldots, \{e_n\}\times\Delta_n$.  Hence these are also facets.

It follows $\tilde\eta_i$ and $\tilde\nu_j$ define facets of $\Delta$,
and these are {\it all\/} of the facets since $\Delta$ is given by
(3.3).  This proves that we have found the generators of the
1-dimensional cones of the fan of $Y$.

We next turn our attention to the coordinate ring of $Y$, which is the
polynomial ring
$$
R = \C[x_1,\ldots,x_{n+r},y_0,\dots,y_n],
$$
where $x_i$ corresponds to $\tilde\eta_i$ and $y_j$ corresponds to
$\tilde\nu_j$.  To determine the grading on $R$, note that the
$\P^n$-fibration $p:Y \to X$ gives an exact sequence
$$
0 \longrightarrow A_{n-1}(X) \mapname{p^*} A_{2n-1}(Y)
\longrightarrow \Z \longrightarrow 0. \leqno{(3.4)}
$$
In terms of $p:Y\to X$, we can think of the $x_i$ as variables coming
from the base and the $y_j$ as variables on the fiber.  To make this
more precise, let the torus invariant divisors on $Y$ corresponding to
$\tilde\eta_i$ and $\tilde\nu_j$ be $\widetilde D_i$ and $\widetilde
D'_j$ respectively.  Then $\widetilde D_i$ is the pullback of the
torus invariant divisor $D_i$ on $X$ corresponding to $\eta_i$, and
$\widetilde D'_j$ induces the hyperplane class on each fiber.  In
particular, $\deg(x_i) = [\widetilde D_i] \mapsto 0$ and $\deg(y_j) =
[\widetilde D'_j] \mapsto 1$ in (3.4).

We next have the following important lemma.

\proclaim Lemma 3.5. For each $j = 0,\ldots,n$, we have ${\cal
O}_Y\big(\widetilde D'_j\big)\otimes p^*({\cal L}_j) \simeq {\cal
O}_Y(1)$.

\noi {\bf Proof.}\ \  The integers $a_{ij}$ in (3.1) mean that ${\cal
L}_j \simeq \ox\big(\sum_{i=1}^{n+r} a_{ij} D_i\big)$ on $X$.  It
follows that on $Y$, we have
$$
{\cal O}_Y\big(\widetilde D'_j\big)\otimes p^*({\cal L}_j) \simeq
{\cal O}_Y\big(\widetilde D'_j + \textstyle{\sum_{i=1}^{n+r}} a_{ij}
\widetilde D_i\big).
$$
When $j = 0$, the polytope corresponding to this divisor is {\it
precisely\/} $\Delta$ by (3.3), which proves the lemma in this case.
If $j > 0$, we have $e_j \in \Z^n\oplus M$, and the divisor of the
corresponding character $\chi^{e_j}$ is
$$
\eqalign{{\rm div}(\chi^{e_j}) &= \textstyle{\sum_{i=1}^{n+r}} \langle
e_j,\tilde\eta_i\rangle \widetilde D_i + \textstyle{\sum_{k=0}^{n}}
\langle e_j,\tilde\nu_k\rangle \widetilde D'_k\cr
&= \textstyle{\sum_{i=1}^{n+r}} (a_{ij}-a_{i0}) \widetilde D_i -
\widetilde D'_0 + \widetilde D'_j\cr
&= \big(\widetilde D'_j + \textstyle{\sum_{i=1}^{n+r}} a_{ij}
\widetilde D_i\big) - \big(\widetilde D'_0 +
\textstyle{\sum_{i=1}^{n+r}} a_{i0} \widetilde D_i\big),\cr}
$$
and the lemma follows immediately.$\ \ \diamond$
\medskip

To see what this lemma says about coordinate rings, let $\alpha_j =
[{\cal L}_j] \in A_{n-1}(X)$ and pick polynomials $F_j \in
S_{\alpha_j}$.  The $F_j$ may may have different degrees in $S$ (since
the $\alpha_j$ need not be equal), but Lemma 3.5 implies that the
polynomials $y_jF_j$ all have the {\it same\/} degree in $R$.  Thus we
can form the single homogeneous polynomial $\sum_{j=0}^n y_jF_j \in R$
which contains all the $F_j$ simultaneously.  This is the essence of
the Cayley trick.

We can now prove the first main result of this section, which gives a
sufficient condition for $\langle F_0,\ldots,F_n\rangle
\subset S$ to have codimension one in the critical degree $\rho$.

\proclaim Theorem 3.6. Let $X$ be a complete simplicial toric variety
of dimension $n$, and assume $F_j \in S_{\alpha_j}$, for $j =
0,\ldots,n$, where $\alpha_j$ is ample and the $F_j$ don't vanish
simultaneously on $X$.  If $\rho = \rho(F)$ is the critical degree of
the $F_j$, then
$$
\dim_\C(S_\rho/\langle F_0,\ldots,F_n\rangle_\rho) = 1.
$$

\noi {\bf Proof.}\ \ If we pick ample line bundles ${\cal L}_j$ on $X$
such that $\alpha_j = [{\cal L}_j] \in A_{n-1}(X)$, then we get the
toric variety $Y = \P({\cal L}_0\oplus\cdots\oplus{\cal L}_n)$
described above.  As remarked after the proof of Lemma 3.5, the
polynomials $y_jF_j$ all have the same degree in the coordinate ring
$R$ of $Y$.  This degree is the ample class $[{\cal O}_Y(1)] \in
A_{2n-1}(Y)$, which we will denote by $\gamma$.  Thus we can define
the homogeneous polynomial
$$
F = y_0 F_0 + \cdots + y_n F_n \in R_\gamma.
$$
Let $W \subset Y$ be the hypersurface defined by $F = 0$.  The idea of
the Cayley trick is that this hypersurface should be closely related
to the complete intersection $F_0 = \cdots = F_n = 0$ on $X$.  Since
the intersection is empty in our situation, we expect $W$ to be
especially simple.

We next check that all of the relevant hypotheses of [BC] are
satisfied.  We know that $\gamma$ is ample, and $Y$ is simplicial
since it is a $\P^n$-bundle over the simplicial toric variety $X$.  To
show that $W$ is quasi-smooth (as defined in Section 3 of [BC]), note
that among the partial derivatives of $F$, we have
$$
{\partial F \over \partial y_j} = F_j.\leqno{(3.7)}
$$
Since the $F_j$ don't vanish simultaneously on $X$, these partials of
$F$ can't vanish simultaneously on $Y$, which proves that $W$ is
quasi-smooth.

The {\it primitive cohomology\/} of $W$ is defined by the exact
sequence
$$
H^{2n-1}(Y) \longrightarrow H^{2n-1}(W) \longrightarrow PH^{2n-1}(W)
\longrightarrow 0
$$
(with coefficients in $\C$).  To prove Theorem 3.6, we will compute
$PH^{2n-1}(W)$ topologically, using $W \hookrightarrow Y \to X$, and
algebraically, using the Jacobian ideal of $F$.

In the composition $W \hookrightarrow Y \to X$, the fiber over a point
of $X$ with coordinates $t_1,\ldots,t_{n+r}$ is the subset of $\P^n$
defined by $\sum_{j=0}^n y_j F_j(t_1,\ldots,t_{n+r}) = 0$.  Since the
$F_j$ don't vanish simultaneously on $X$, it follows that the fiber is
a hyperplane $\P^{n-1} \subset \P^n$.  Topologically, this means we
have a map of fibrations
$$
\matrix{\P^{n-1} & \hookrightarrow & \P^n\cr
	\downarrow && \downarrow \cr
	W & \hookrightarrow & Y\cr
	\downarrow && \downarrow\cr
	X & = & X.\cr}
$$
For each fibration, we get the usual spectral sequence, and the map
between the spectral sequences is surjective at $E_2$ because
$H^q(\P^n) \to H^q(\P^{n-1})$ is surjective for all $q$.  It follows
that $H^{2n-1}(Y) \to H^{2n-1}(W)$ is surjective, so that
$PH^{2n-1}(W)$ vanishes.

We can also compute the Hodge components of $PH^{2n-1}(W)$ using [BC].
In particular, the exact sequence from Theorem 10.13 of [BC] gives an
exact sequence
$$
0 \to H^{2n-2}(Y) \to H^{2n}(Y) \to (R/J(F))_{(n+1)\gamma-\tilde\beta}
\to PH^{n-1,n}(W)\, (= 0) \to 0,\leqno{(3.8)}
$$
where $J(F) = \langle \partial F/\partial x_i,\partial F/\partial
y_j\rangle$ is the Jacobian ideal of $F$ and $\tilde\beta =
\sum_{i=1}^{n+r} \deg(x_i) + \sum_{j=0}^n \deg(y_j)$.  However,
$\gamma = \deg(y_j) + \alpha_j$ for all $j$ by Lemma 3.5, so that
$$
\eqalign{
(n+1)\gamma-\tilde\beta &= \textstyle{\sum_{j=0}^n} (\deg(y_j) +
	\alpha_j) - \textstyle{\sum_{i=1}^{n+r}} \deg(x_i) -
	\textstyle{\sum_{j=0}^n} \deg(y_j)\cr
	&= \textstyle{\sum_{j=0}^n} \alpha_j -
	\textstyle{\sum_{i=1}^{n+r}} \deg(x_i) = \rho.\cr}
$$
In the map $A_{2n-1}(Y) \to \Z$ of (3.4), we know that $\rho \mapsto
0$ and $\deg(y_j) \mapsto 1$.  This implies $R_\rho = S_\rho$.
Furthermore, by (3.7), the Jacobian ideal is $J(F) = \langle F_j,
\partial F/\partial x_i\rangle$, and $J(F)_\rho = \langle
F_0,\ldots,F_n\rangle_\rho$ follows since $\partial F/\partial x_i =
\sum_{j=0}^n y_j\partial F_j/\partial x_i$.  Then (3.8) tells us that
$$
\dim_\C(S_\rho/\langle F_0,\ldots,F_n\rangle_\rho) = h^{2n}(Y) -
h^{2n-2}(Y).
$$
However, since the spectral sequence for the fibration $\P^n \to Y \to
X$ degenerates at $E_2$ (both base and fiber have cohomology only in
even degrees), we see that if $q \le n$, then
$$
h^{2q}(Y) = \textstyle{\sum_{k=0}^q} h^{2k}(X) h^{2(q-k)}(\P^n) =
h^0(X) +  h^2(X) + \cdots + h^{2q}(X).
$$
This easily implies $\ \dim_\C(S_\rho/\langle F_0,\ldots,F_n\rangle_\rho) =
h^{2n}(X) = 1$, and the
theorem is proved.$\ \ \diamond$
\bigskip

We can now prove Theorem 0.3 from the Introduction.  The first part of
the theorem claims that the toric residue map
$$
\res_F : S_\rho/\langle F_0,\ldots,F_n\rangle_\rho \longrightarrow \C
\leqno{(3.9)}
$$
is an isomorphism.  Since $X$ is simplicial, every $n$-dimensional
$\sigma \in \Sigma$ is simplicial, so that by Theorem 0.2, we have
$\Delta_\sigma \in S_\rho$ such that $\res_F(\Delta_\sigma) = \pm1$.
Then Theorem 3.6 immediately implies that (3.9) is an isomorphism.

Turning to the second part of Theorem 0.3, we need to show that
$$
x_i \cdot S_\rho \subset \langle F_0,\ldots,F_n\rangle,\quad i =
1,\ldots,n+r. \leqno{(3.10)}
$$
To prove this, let $\sigma$ be a $n$-dimensional cone of $\Sigma$
containing $\eta_i$, and renumbering as in \S2, we can assume that $i
\le n$.  Then Cramer's Rule, applied to the equations (2.2),
shows that $x_i\cdot \Delta_\sigma \subset \langle
F_0,\ldots,F_n\rangle$.  But the previous paragraph implies $S_\rho =
\C\cdot \Delta_\sigma + \langle F_0,\ldots,F_n\rangle_\rho$, and then
(3.10) follows immediately.  This completes the proof of Theorem 0.3.
\bigskip

\noindent {\bf Remark 3.11.} As a consequence of these results, we can
describe an algorithm for computing toric residues when $X$ is
complete and simplicial and $F_j \in S_{\alpha_j}$ for $\alpha_j$
ample.  First, pick a Gr\"obner basis for $\langle
F_0,\ldots,F_n\rangle$ (using a convenient monomial order on $S$).
Given a polynomial $H \in S$, we can then compute its normal form,
denoted ${\it normalform}(H)$.  Since $\langle
F_0,\ldots,F_n\rangle_\rho \subset S_\rho$ has codimension 1, an easy
argument shows that the normal forms of elements of $S_\rho$ are
multiples of the monomial $x^\alpha$ which is the least (relative to
the chosen monomial order) among the monomials of degree $\rho$ not in
$\langle F_0,\ldots,F_n\rangle_\rho$.

Then choose a $n$-dimensional cone $\sigma$, say with generators
$\eta_{i_j}$, and pick a basis $m_i$ of $M$ such that $\det(\langle
m_i,\eta_{i_j}\rangle) > 0$.  If we use this basis to construct the
Euler form $\Omega$, then by the remarks made at the end of \S2, the
determinant $\Delta_\sigma \in S_\rho$ satisfies
$\res_F(\Delta_\sigma) = 1$.  Finally, let $c_\sigma$ be the nonzero
constant such that ${\it normalform}(\Delta_\sigma) = c_\sigma
x^\alpha$.

Given these ``preprocessing'' steps, we can now describe the
algorithm: given $H\in S_\rho$, its toric residue is given by the
quotient
$$
\res_F(H) = {c \over c_\sigma},
$$
where ${\it normalform}(H) = c\, x^\alpha$.  This follows because $H
\equiv c\, x^\alpha \bmod \langle F_0,\ldots,F_n\rangle$ and
$\Delta_\sigma \equiv c_\sigma\,x^{\alpha} \bmod \langle
F_0,\ldots,F_n\rangle$ imply $H \equiv (c/c_\sigma)\, \Delta_\sigma
\bmod \langle F_0,\ldots,F_n\rangle$.
\bigskip

In the final part of this section, we will discuss the hypotheses of
the Codimension One Theorem (Theorem 3.6) and the Residue Isomorphism
Theorem (Theorem 0.3).  In proving both of these results, we assumed
that the degrees of $F_0,\ldots,F_n$ were ample classes in
$A_{n-1}(X)$ (this was needed in order to use the results of [BC]).
We suspect that these theorems should hold under the weaker hypothesis
that the degrees are $\Q$-ample.  In fact, there is an even weaker
hypothesis which leads to the following conjecture generalizing the
Codimension One Theorem.

\proclaim Conjecture 3.12. If $X$ is a complete simplicial toric
variety and $F_0,\ldots,F_n \in B(\Sigma)$ are homogeneous polynomials
which don't vanish simultaneously on $X$, then
$$
\dim_\C(S_\rho/\langle F_0,\ldots,F_n\rangle_\rho) = 1,
$$
where as usual $\rho$ is the critical degree of $F_0,\ldots,F_n$.

Recall from \S2 that $B(\Sigma)$ is the ideal generated by the
monomials $\hat x_{\sigma} = \Pi_{\eta_i\notin\sigma} x_i$ for all
$\sigma\in\Sigma$ and that $S_{\alpha_j} \subset B(\Sigma)$ when
$\alpha_j$ is ${\bf Q}$-ample (see (2.1)).  Thus Theorem 3.6 is a
special case of Conjecture 3.12.

One useful observation is that Conjecture 3.12 implies the conclusions
of the Residue Isomorphism Theorem remain true.

\proclaim Proposition 3.13.  Let $X$ be a complete simplicial toric
variety, and let $F_0,\ldots,F_n \in B(\Sigma)$ be homogeneous
polynomials which don't vanish simultaneously on $X$.  If
Conjecture~3.12 is true for $X$ (i.e., if $\dim_\C (S_\rho/\langle
F_0,\ldots,F_n\rangle_\rho) = 1$), then:
\vskip0pt
\item{(i)} The toric residue map $\res_F : S_\rho/\langle
F_0,\ldots,F_n\rangle_\rho \to \C$ is an isomorphism.
\item{(ii)} For each variable $x_i$, $0 \le i \le n+r$, we have
$x_i\cdot S_\rho \subset \langle F_0,\ldots,F_n\rangle$.
\vskip0pt

\noi {\bf Proof.}\ \ The argument is identical to what we used to
derive Theorem 0.3 from Theorem~3.6.  This is because, as we observed
in Remark 2.12 (ii), Theorem 0.2 still applies under the assumption
$F_i \in B(\Sigma)$.$\ \ \diamond$
\medskip

As evidence for Conjecture 3.12, we present the following examples.
\bigskip

\noindent {\bf Examples 3.14.}\ (i)\ If $X = \P(q_0,\ldots,q_n)$ is a
weighted projective space with coordinate ring $S =
\C[x_0,\ldots,x_n]$, then $B(\Sigma)$ is the ideal $\langle
x_0,\ldots,x_n\rangle$, so that $F_i \in B(\Sigma)$ means that $F_i$
has positive degree.  Hence Conjecture 3.12 follows easily by standard
commutative algebra because $F_0,\ldots,F_n$ form a regular sequence
in $S$ (since they don't vanish simultaneously on $X$).

	For a specific example, consider $X = \P(1,1,1,1,3,3,5)$ and
suppose that $F_0,\ldots,F_6$ have degrees 3, 6, 6, 6, 6, 5, 4
respectively.  The critical degree is $\rho = 21$, so that
$S_{21}/\langle F_0,\ldots,F_6\rangle_{21} \simeq \C$ in this case.
Note that the $F_i$ are {\it not\/} Cartier, though they are certainly
$\Q$-ample.  This is an example from mirror symmetry which arises in
connection with certain $(0,2)$ string theories---see [DK] for more
details.
\medskip

\noi (ii)\ For another example where the degrees of the $F_i$ are
$\Q$-ample but not Cartier, consider the toric surface $X$
corresponding to the fan in $\R^2$ determined by the vectors
$$
\eta_1 = (1,0);\quad \eta_2 = (0,1);\quad \eta_3 =
(-1,1);\quad \eta_4 = (-1,-1);\quad \eta_5 = (1,-1).
$$
Note that $X$ is singular since $\eta_3$, $\eta_4$ and $\eta_4$,
$\eta_5$ don't span all of $\Z^2$, though $X$ is certainly simplicial.

If we let the variables $x,y,z,t,u$ correspond to
$\eta_1,\ldots,\eta_5$, then the exceptional set $Z \subset \C^5$ is
defined by the ideal $B(\Sigma) = \langle ztu , xtu , xyu , xyz , yzt
\rangle$, that is,
$$
Z= \{x=z=0 \} \cup \{z=u=0 \} \cup\{y=u=0 \} \cup\{x=t=0 \} \cup
\{y=t=0 \}.
$$
Thus, $X \simeq (\C^5 - Z) / (\C^*)^3$.  Furthermore, one can show
that $A_1(X) \simeq \Z^3$ and that we get a grading in the polynomial
ring $S = \C[x,y,z,t,u]$ with
$$
\displaylines{
{\rm deg}(x) = (1,1,-1);\quad
{\rm deg}(y) = (-1,1,1);\quad
{\rm deg}(z) = (1,0,0);\cr
{\rm deg}(t) = (0,1,0);\quad
{\rm deg}(u) = (0,0,1).\cr}
$$
Thus, the sum of the degrees of the variables is $\beta = (1,3,1)$.

We next characterize ample divisors on $X$.  First, one checks
that a class $(a,b,c) \in \Z^3 \simeq A_1(X)$ lies in ${\rm Pic}(X)
\subset A_1(X)$ (i.e., the divisor $a\,D_3 + b\,D_4 + c\,D_5$ is
Cartier) if and only if $a \equiv b \equiv c \bmod 2$.  Then it is
straightforward to verify (using [F, \S3.3-4]) that a Cartier class
$(a,b,c)$ is ample if and only if
$$
b>a>0\quad\quad{\rm and}\quad\quad b>c>0.\leqno{(3.15)}
$$
For an arbitrary $(a,b,c)$, these inequalities tell us when the
corresponding class is $\Q$-ample.

Now consider the polynomials
$$
F_0 = xy^2z^3;\quad F_1 = x^2yu^3 + yz^2t^2u + xt^2u^3 + y^2z^3t;\quad
F_2 = zt^3u^2 + xt^2u^3 + y^2z^3t.
$$
They are homogeneous and $\deg F_0 = (2,3,1)$, $\deg F_1 = \deg F_2 =
(1,3,2)$.  One can check that the common zeros of $F_0$, $F_1$ and
$F_2$ in $\C^5$ are contained in the set $Z$ and therefore the
corresponding divisors on $X$ have empty intersection.  None of these
divisors are Cartier, but they are clearly $\Q$-ample by (3.15), and
their critical degree is given by:
$$
\rho = (2,3,1) + (1,3,2) + (1,3,2) - (1,3,1) = (3,6,4).
$$
There are 22 monomials of degree $\rho$, and computing the normalform
of each monomial (as in Remark 3.11), we find that the normalforms are
all multiples of the same monomial (for example, if we use graded
reverse lex with $x > y > z > t > u$, the normalforms are all
multiples of $x^3t^3u^7$).  Thus $\langle F_0,F_1,F_2 \rangle_{\rho}$
has codimension one in $S_{\rho}$.
\medskip

\noi (iii) We next give an example where $F_i \in B(\Sigma)$ for all
$i$ but their degrees are not $\Q$-ample classes.  We use the
same singular toric surface $X$ as in (ii), but this time we consider
the polynomials
$$
F_0 = ztu;\quad F_1 = yzt + xyu;\quad F_2 = xyz + xtu.
$$
These are homogeneous with degrees $\deg F_0 = (1,1,1)$, $\deg F_1 =
(0,2,1)$, and $\deg F_2 = (1,2,0)$.  One can check that $F_0$, $F_1$
and $F_2$ don't vanish simultaneously on $X$, and by the ampleness
criterion (3.15), none of their degrees are $\Q$-ample, although
$F_0,F_1,F_2$ all lie in $B(\Sigma)$.  The critical degree is
$$
\rho = (1,1,1) + (0,2,1) + (1,2,0) - (1,3,1) = (1,2,1).
$$
Computing normalforms of the four monomials of degree $\rho$ reveals
that $\langle F_0,F_1,F_2 \rangle_{\rho}$ has codimension one in
$S_{\rho}$.
\medskip

\noi (iv)\ Finally, we give an example to show what can go wrong if
not all of the $F_i$ are contained in $B(\Sigma)$.  Let $X= \P^1\times
\P^1$.  Here, it is well known that the homogeneous coordinate ring of
$X$ is $S = \C[x,y,z,t]$, with the usual bigrading
$$
{\rm deg}(x) ={\rm deg}(y) = (1,0);\quad
{\rm deg}(z) ={\rm deg}(t) = (0,1).
$$
Also, $B(\Sigma) = \langle xz,xt,yz,yt\rangle$.

We now let
$$
F_0 = (x+y)^2;\quad F_1 = xz;\quad F_3 = yt.
$$
Thus $\deg(F_0) = (2,0)$ and $\deg(F_1) = \deg(F_2) = (1,1)$.  It is
easy to check that $F_0, F_1, F_2$ don't vanish simultaneously on $X$.
Moreover, the divisors defined by $F_1, F_2$ are ample (a polynomial
of degree $(a,b)$ defines an ample divisor if and only if $a>0$ and
$b>0$), while $F_0 \notin B(\Sigma)$.

The critical degree in this case is $\rho = (2,0)$ since the sum of
the degrees of the variables is $\beta = (2,2)$.  There are three
monomials of degree $(2,0)$: $x^2$, $y^2$ and $xy$, and any two of them
are linearly independent modulo the ideal $\langle F_0,F_1,F_2
\rangle$.  Thus $\langle F_0,F_1,F_2\rangle_\rho$ does {\it not\/}
have codimension one in $S_\rho$.  Note also that {\it no\/} monomial
of degree $(3,0)$ is in the ideal, which shows that $x\cdot S_\rho
\not\subset \langle F_0,F_1,F_2\rangle$.  Hence the second part of
Proposition 3.13 fails as well.
\bigskip

\noi {\bf Remarks 3.16.}\ (i)\ Notice that if the $F_i$ don't all lie
in $B(\Sigma)$, then we can no longer express the $F_i$ as in (2.2),
so that the definition of $\Delta_{\sigma}$ makes no sense.  Thus,
even if $\langle F_0,\ldots,F_n\rangle_\rho$ has codimension one in
$S_\rho$, the second part of Proposition 3.13 could fail.  For an
example of how this can happen, consider the toric variety $X$ of
Example 3.14 (ii), this time using the polynomials
$$
F_0 = ztu;\quad F_1 = yzt + xyu;\quad F_2 = xyz + xtu + zt^2.
$$
These are very similar to what we used in Example 3.14 (iii)---the
only difference is that $F_2$ has an extra $zt^2$ term.  As in that
example, $\deg F_0 = (1,1,1)$, $\deg F_1 = (0,2,1)$, and $\deg F_2 =
(1,2,0)$, and they don't vanish simultaneously on $X$.  Note also that
$F_0, F_1 \in B(\Sigma)$ but $F_2\not\in B(\Sigma)$ because of the
$zt^2$ term.  The critical degree is still $(1,2,1)$, and an easy
computation shows that $\langle F_0,F_1,F_2 \rangle_{\rho}$ still has
codimension one in $S_{\rho}$.  However, in this case, one can also
compute that
$$
x\cdot xyzu \notin \langle F_0,F_1,F_2 \rangle.
$$
Since $xyzu \in S_\rho$, we have $x\cdot S_\rho \not\subset \langle
F_0,F_1,F_2 \rangle$, so that the second part of Proposition 3.13
fails in this case.
\medskip

\noi (ii)\ One question we have not investigated is whether the
simplicial hypothesis is needed in Conjecture 3.12 and Proposition
3.13.  For example, if $X$ is an arbitrary complete toric variety,
then Conjecture~3.12 and the first part of Proposition 3.13 are true
when the degrees of the $F_i$ are the same ample class---this is Theorem 5.1
of [C2].
\bigskip
\bigskip

\noi {\bf \S4. Global Residues as Sums of Local Residues}
\bigskip

In this section we will show that for simplicial toric varieties, the
toric residue may be computed as a sum of local Grothendieck residues.
The toric setting is not essential here and, in fact, it is convenient
to work with the more general notion of a $V$-manifold or orbifold
(see [B], [Sa]).  We begin with a review of the theory of residual
currents.
\medskip

{\bf Residual Currents on $V$-Manifolds. } We recall that, by results
of Prill [P], if an $n$-dimensional complex variety $X$ is a
$V$-manifold, then for every $x \in X$ there exists a finite subgroup
$G\subset GL(n,\C)$ such that for some neighborhood $W$ of $x \in X$,
we have $(W,x) \simeq (U/G,0)$, where $U$ is a $G$-invariant
neighborhood of $0 \in C^n$.  Furthermore, $G$ is {\it small\/} (no
$g\in G$ has $1$ as an eigenvalue of multiplicity $n-1$) and is unique
up to conjugacy.  Such a local presentation $(W,x) \simeq (U/G,0)$ is
called a {\it standard model\/}.

A simplicial toric variety $X$ is an example of a $V$-manifold.
Indeed, with the notation of \S 2, we may cover $X$ with open sets
$X_{\sigma} \simeq \C^n/G(\sigma)$ and it is easy to verify that
$G(\sigma)$ is a small subgroup (see [BC, 3.5]).

It is shown in [St, 1.8] that if $X$ is a $V$-manifold and $(W,x)
\simeq (U/G,0)$ is a standard model, then
$$
\Gamma(W,\widehat\Omega^p_X) \simeq
\Gamma(U,\Omega_{\C^n}^p)^G\leqno{(4.1)}
$$
where, as before, $\widehat\Omega^p_X$ denotes the sheaf of Zariski
$p$-forms on $X$, and the superscript $G$ indicates the subspace of
$G$-invariant forms.  Similarly (see [Sa] and [B]), we consider the
sheaves $\EE^{p,q}_X$ of $\cinf$ forms on $X$ of bidegree $(p,q)$.
They are associated with the presheaves which assign to an open set
$W\subset X$, which is part of a standard model $(W,x)\simeq (U/G,0)$,
the group $\EE^{p,q}_X(W) = \Gamma(U,\EE^{p,q})^G$, where $\EE^{p,q}$
is the sheaf of $\cinf$ $(p,q)$-forms on $\C^n$.  The restriction maps
for these presheaves are defined as follows: if $(W',x) \simeq
(U'/G',0)$ is another standard model and $W'\subset W$ then by [P,
Theorem 2], there exists a linear map $h\in GL(n,\C)$ such that
$h(U')\subset U$ and $G' = h^{-1} G h$.  We then set $r^W_{W'} = h^* :
\EE^{p,q}_X(W) \to \EE^{p,q}_X(W')$.  Note also that any
element in $GL(n,\C)$ commutes with the differential operator
$\bar\partial$ acting on $\Gamma(U,\EE^{p,q})$, which means that we
can define an operator $\bar\partial : \EE^{p,q}_X \to \EE^{p,q+1}_X$.

We denote by $\Gamma_c(W,\EE^{p,q}_X)$ the space of sections of
$\EE^{p,q}_X$ with compact support in $W$.  For a standard model
$(W,x) \simeq (U/G,0)$, we have $\Gamma_c(W,\EE^{p,q}_X) \simeq
\Gamma_c(U,\EE^{p,q})^G$, and it carries a natural Fr\'echet topology
as a subspace of $\Gamma_c(U,\EE^{p,q})$.  We will denote by $\curr
pqX$ the sheaf of $(p,q)$-{\it currents\/} on $X$, i.e., the sheaf
which associates to any open set $W$ of $X$, the space $\curr pqX(W)$
of continuous linear funcionals on $\Gamma_c(W,\EE^{n-p,n-q}_X)$.

\proclaim Lemma 4.2.   If $(W,x) \simeq (U/G,0)$ is a standard model,
then
$$
\curr pqX(W) \simeq \curr pq{}(U)^G,
$$
where $\curr pq{}$ is the sheaf of $(p,q)$-currents on $\C^n$ and the
action of $G$ on $\curr pq{}(U)$ is the natural one:
$$
(gT)(\alpha) = T(g^*\alpha),\quad T\in \curr pq{}(U)\ \ \hbox{and}\ \
\alpha\in \Gamma_c(U,\EE^{n-p,n-q}).
$$

\noi{\bf Proof.}\ \ The space $\curr pqX(W)$ is by definition
the continuous dual of $\Gamma_c(W,\EE^{n-p,n-q}_X) =
\Gamma_c(U,\EE^{n-p,n-q})^G$, and any continuous linear map
$$
T: \Gamma_c(U,\EE^{n-p,n-q})^G \longrightarrow \C
$$
extends to a $G$-invariant continuous linear map defined on all of
$\Gamma_c(U,\EE^{n-p,n-q})$ by the formula
$$
T(\alpha) = T(\alpha^G),\quad  {\rm where}\ \alpha^G = \frac {1}{|G|}
\sum_{g\in G} g^*\alpha.\leqno{(4.3)}
$$
Conversely, every $G$-invariant linear functional $T$ on
$\Gamma_c(U,\EE^{n-p,n-q})$ satisfies
$$
T(\alpha) = \frac {1}{|G|}\sum_{g\in G} (gT)(\alpha) =
\frac {1}{|G|}\sum_{g\in G} T(g^*\alpha) = T(\alpha^G)\leqno{(4.4)}
$$
and, thus, is in the image of (4.3).   $\ \ \diamond$
\medskip

It will be convenient to assign to $T \in \curr pq{}(U)^G$ the element
$(1/|G|)\,T \in \curr pqX(W)$. With this convention, the $G$-invariant
continous linear operator defined by integration $\int_U :
\Gamma_c(U,\EE^{n,n}) \to \C$ gives rise, when $(W,x)\simeq (U/G,0)$
is a standard model, to the usual definition of integration for
sections $\alpha\in\Gamma_c(W,\EE^{n,n}_X)$:
$$
\int_W\alpha = {\frac 1{|G|}}\ \int_U \alpha.
$$
It is clear that this definition is independent of the choice of
standard model.  Moreover, the existence of $\cinf$ partitions of unity
on $X$ (see [B]) implies that we can define the integral for compactly
supported sections of $\EE^{n,n}_X$ over any open set of $X$.

Similarly, given a $G$-invariant form $\alpha\in
\Gamma(U,\EE^{p,q})^G$, integration against $\alpha$ defines a
$G$-invariant current $I(\alpha)\in \curr pq{}(U)^G$.  Thus
$$
I(\alpha)(\beta) = \int_U \alpha\wedge\beta\quad \hbox{for}\quad
\beta\in\Gamma_c(U,\EE^{n-p,n-q}).
$$
The corresponding current in $\curr pqX(W)$ will also be denoted by
$I(\alpha)$ and we have $I(\alpha)(\beta) = (1/|G|)\int_U
\alpha\wedge\beta = \int_W \alpha\wedge\beta$.

We extend the definition of $\bar\partial$ to the space of currents by
the formula:
$$
(\bar\partial T)(\beta) := (-1)^{p+q}\ T(\bar\partial \beta),\quad
T\in \curr pqX(W)\ \ \hbox{and}\ \ \beta \in
\Gamma_c(W,\EE_X^{n-p,n-q-1}).
$$

\proclaim Proposition 4.5.  Let $X$ be a compact, connected
$V$-manifold.  Then:
\vskip0pt
\item{(i)}\ The diagram
$$\matrix{0&\longrightarrow&
\widehat\Omega_X^p&
\longrightarrow&
\EE^{p,0}_X&
\mapname{\dbar}&
\cdots&
\mapname{\dbar}&
\EE^{p,n}_X&
\longrightarrow&
0\cr
&&\Vert
&&\downarrow\, I
&&&&\downarrow\, I
&&\cr
0&\longrightarrow&
\widehat\Omega_X^p&
\mapname{I}&
\curr p0X&
\mapname{\dbar}&
\cdots&
\mapname{\dbar}&
\curr pnX&
\longrightarrow&0\cr}
$$
commutes and its rows are exact.
\item{(ii)}\ The following diagram commutes and all maps are isomorphisms:
$$
\matrix{
&&H^{n,n}_{\dbar}(X)&&\cr
&\noe {\eta}&&\hskip-.1truein\soe {\int_X}&\cr
H^n(X,\widehat\Omega^n_X)&&\mapdown {I}&&\hskip-.1truein\C\cr
&\soe {{}'\!\eta}&&\hskip-.1truein\noe{{\rm Ev}_1}&\cr
&&H^{n}_{\dbar}(\Gamma(X,\curr n\cdot{X}))&&\cr}
$$

\noi{\bf Proof.}\ \ The commutativity of the first diagram is a
consequence of the sign convention in the definition of
$\bar\partial$.  Exactness follows from the corresponding statements
in the smooth case.  We illustrate this for the bottom row.

Let $\alpha\in\widehat\Omega^p_{X,x}$ be such that $I(\alpha)=0$.  We
represent $\alpha$ by a $G$-invariant holomorphic $p$-form $\tilde
\alpha$ on $U$ where $(U/G,0)\simeq (W,x)$ is a standard model. By
Lemma 4.2, $I(\tilde\alpha) = 0$ as an element in $\curr p0{}(U)^G
\subset \curr p0{}(U)$. Consequently, by exactness in the smooth case,
$\tilde \alpha = 0$ and, a fortiori, $\alpha=0$.

Suppose now that $T\in \curr pq{X,x}$ is $\bar\partial$-closed.
Again, we represent $T$ by a $G$-invariant current $\tilde T \in \curr
pq{}(U)^G$ satisfying $\bar\partial\tilde T=0$.  We may replace $U$ by
a smaller $G$-invariant neighborhood $U'$ of $0\in\C^n$ where $\tilde
T = \bar\partial \tilde S$, $\tilde S \in \curr p{q-1}{}(U')$.  As in
(4.4), since $\bar\partial$ is a $G$-invariant operator, $\tilde T =
\bar\partial \tilde S^G\,,$ where $\tilde S^G$ is the $G$-invariant
current $\tilde S^G = (1/|G|) \sum_{g\in G} g\tilde S$.  Thus $T =
\bar\partial S^G$, for the induced element in $\curr p{q-1}{X,x}$.

To prove (ii) we note that the sheaves $\EE^{p,q}_X$ and $\curr pqX$
are fine and, consequently, the rows in the diagram in (i) give fine
resolutions of the sheaf $\widehat\Omega^p_X$.  Now, taking $p=n$, the
usual proof of Dolbeault's Theorem gives the isomorphisms $\eta$ and
${}'\!\eta$.  The isomorphism $I$ is deduced from the map at the level
of sheaves and the commutativity follows from (i). Clearly $I$ maps
the cohomology class of a $\dbar$-closed $(n,n)$-form $\alpha$ to the
$(n,n)$-current defined by integration of compactly supported $\cinf$
functions against $\alpha$.

Stokes' Theorem for $V$-manifolds [B] implies that integration over
$X$ defines an isomorphism $\int_X\colon H^{n,n}_{\dbar}(X) \to \C$
and the map ${\rm Ev}_1 : H^{n}_{\dbar}(\Gamma(X,\curr n\cdot{X})) \to
\C$ is defined by evaluation of a (global) current on the constant
function $1_X$.  The commutativity of the right triangle then follows
from the relation $\int_X \alpha = I([\alpha])(1_X)$, $[\alpha]\in
H^{n,n}_{\dbar}(X)$.  $\ \ \diamond$
\medskip

We now bring into the picture Multiple Residue and Principal Value
currents (as in [CH] and [Di]).  Let $D_1,\ldots, D_k$ be reduced Weil
divisors on the $V$-manifold $X$.  For each $D_j$, some multiple is a
Cartier divisor (since $X$ is a $V$-manifold), so that $D_j$ may be
given locally as the support of the zero set of a holomorphic
function.  Let $\omega$ be a semimeromorphic $(p,q)$-form on $X$ with
poles on $D = \cup^k_{i=1} D_i$.  This means that $\omega$ can be
locally written as $\omega'/f$ with $\omega'$ a $\cinf$ $(p,q)$-form
and $f$ a holomorphic function such that $\{f = 0\} \subset D$.

Suppose for a moment that $X$ is smooth and we are given (not
necessarily minimal) equations $f \sb i \in \Gamma (U, {\cal O} \sb
{X})$ for each hypersurface $D_i$, $i = 1,\ldots,k$, on some open
subset $U$. For any $\cinf$ form $\alpha$ (resp.~$\beta$) with compact
support contained in $U$ and bidegree $(n-p,n-q-k)$
(resp.~$(n-p,n-q-(k-1))$), we define:
$$
R_D [\omega] (\alpha) =
R_{D_1, \ldots, D_k} [\omega] (\alpha)
= \lim \sb {\delta
\to 0} \int \sb {T \sb \delta (f)} \omega \wedge \alpha
$$
and
$$
RP_D [\omega] (\beta) = R_{D_1,\ldots,D_{k-1}}
P_{D_k} [\omega] (\beta) = \lim \sb
{\delta \to 0} \int \sb {D \sb \delta (f)} \omega \wedge \beta,
$$
where
$$
\displaylines{ T \sb {\delta} (f) = \{ x \in U : |f \sb i (x) | =
\epsilon_i (\delta)\,, \ \ 1 \le i \le k \}\cr
 D \sb \delta (f) = \{ x \in U :
| f \sb i (x) | = \epsilon_i (\delta)\,, \ \ 1 \le i \le k-1, \ \
|f \sb k (x) | > \epsilon_k (\delta)\}\cr}
$$
are conveniently oriented semianalytic tubes and the $k$ functions
$\epsilon_i : (0,1) \to {\bf R} \sb {+}$ are analytic and satisfy
$\lim \sb {\delta \to 0} ({\epsilon_j (\delta)}/ {\epsilon^q \sb {j+1}
(\delta)}) = 0$ for all $1 \le j \le k-1$ and all positive integers
$q$.  We call $(\epsilon_1, \ldots,\epsilon_k)$ an {\it admissible
path\/}.

In [CH], Coleff and Herrera show that the above limits exist for any
$\alpha, \beta$.  Moreover, these limits are independent of the
admissible path and the particular equations $f \sb 1, \ldots, f \sb
k$.  Thus, on $U$, we get the {\it multiple residue\/} current $R_D
[\omega]$ of bidegree $(p, q+k)$ and the {\it principal
value\/} current $RP_D[\omega]$ of bidegree $(p, q+k-1)$.  By means of a
${\cinf}$ partition of unity, these local definitions can be collected
to obtain global currents on $X$, also denoted $R_D[\omega]$ and
$RP_D[\omega]$, whose supports verify
$$
\supp(R \sb {D}
[\omega]) \subset  (\cap^{k}_{i=1} D_i)\cap\supp(\omega),\quad
\supp(RP
\sb {D} [\omega]) \subset  ({\cap}^{k-1}_{i=1} D \sb
i)\cap\supp(\omega).
$$

Suppose now that $X$ is a $V$-manifold, $D_1,\ldots,D_k$ reduced Weil
divisors as above, and $W\simeq U/G$ is a standard model.  We denote
by $\tilde D_1,\ldots,\tilde D_k$ the lifted hypersurfaces in $U$.
For any $G$-invariant, semimeromorphic form $\tilde \omega$ on $U$,
with polar set contained in $\tilde D = \cup_{i=1}^k \tilde D_i$, the
currents $R_{\tilde D}[\tilde \omega]$ and $RP_{\tilde D}[\tilde
\omega]$ are also $G$-invariant.  Thus, given a semimeromorphic form
$\omega$ on $W$, we denote by $\tilde \omega$ its lifting to $U$ and
then define:
$$
R_D[\omega] = \frac 1{|G|} R_{\tilde D}[\tilde \omega]
\quad \hbox{and}\quad
RP_D[\omega] = \frac 1{|G|} RP_{\tilde D}[\tilde \omega].
$$
These definitions may again be globalized using a partition of unity
on $X$.  The definition of $RP_D$ and support property stated above
imply that:
$$
RP_D[\omega]\lower.04truein\hbox{$\big\vert$}
\lower.07truein\hbox{$\scriptstyle{X-D_j}$} = 0\quad \hbox{for } j<k,\
\hbox{and}\quad RP_D[\omega] \lower.04truein\hbox{$\big\vert$}
\lower.07truein\hbox{$\scriptstyle{X-D_k}$} = R_{D_1,\ldots,
D_{k-1}}[\omega]\lower.04truein\hbox{$\big\vert$}
\lower.07truein\hbox{$\scriptstyle{X-D_k}$}.\leqno{(4.6)}
$$

The mappings $R_D$ and $RP_D$ associating to any germ of meromorphic
$p$-form $\omega$ with poles contained in $D$, the germ of the
residual currents $R_D[\omega] $ and $RP_D[\omega]$, define sheaf
morphisms making the following diagram commutative:
$$
\matrix{
\widehat\Omega^p_X(*D)&\hskip-.1truein\hfl{RP_D}{}&\curr p{k-1}X\cr
{}\qquad\hskip.2truein{\scriptstyle R_D}\searrow&&\hskip-.6truein\swarrow
{\scriptstyle {\dbar}}\cr {}&\hskip-.1truein\curr pkX&\cr}\leqno{(4.7)}
$$
In particular, $\dbar R_D [\omega] = 0$ for every meromorphic form
$\omega \in \widehat\Omega^p_X(*D)$.

We conclude our discussion of residual currents by defining the local
Grothendieck residue at a point $x$ on a $V$-manifold $X$.  Let
$(W_1,x)\simeq (U_1/G,0)$ be a standard model, and let $W$ be a
relatively compact neighborhood of $x$ such that $W \subset \overline
W \subset W_1$.  Finally, let $U$ be a $G$-invariant neighborhood of
$0$ such that $W\simeq U/G$ and $U\subset\overline U\subset U_1$.
Suppose that $f_1,\ldots,f_n \in \OO(\overline W)$ have $x$ as their
only common zero in $\overline W$.  Pulling-back to $\overline U$, it
follows that the hypersurfaces $\tilde D_i = \{\tilde f_i = 0\}$
intersect only at $0$.  Given now a meromorphic $n$-form $\omega$ on
$W$ with polar set contained in $\dsize \cup^n_{i=1} \{f_i = 0\}$, we
denote by $\tilde\omega$ its pull-back to $U$ and define
$$
\res \sb
x (\omega) =\frac 1{|G|}\res \sb
0 (\tilde \omega), \leqno{(4.8)}
$$
where, as we recall from [GH] for example, the local Grothendieck
residue $\res \sb 0 (\tilde \omega)$ is defined as
$$
\res \sb 0 (\tilde \omega) = \left(\frac 1{2\pi i} \right)^n
\int_{\{z \in U : |\tilde f
\sb i (z)| =\epsilon_i,\,  1 \le i \le n\}} \tilde\omega.
$$
Here $\epsilon_i > 0$ must be chosen so that $\{w \in {\bf C}^n : |w
\sb i| \sb {1 \le i \le n} = \epsilon_i \}$ is contained in the open
set $\tilde f(U)$, $\tilde f = (\tilde f_1,\ldots, \tilde f_n)$, and
$\{z \in \overline U : |\tilde f \sb i (z)| \sb {1 \le i \le n} =
\epsilon_i \} \cap \partial U = \emptyset$. Note that the tube $\{z
\in U : |\tilde f \sb i (z)| = \epsilon_i, \;1 \le i \le n \}$ is
compact, of real dimension $n$, and we orient it with the form $d
({\rm arg} \;\tilde f \sb 1) \wedge \ldots \wedge d({\rm arg} \;\tilde
f \sb n)$.

If $\varphi$ is a $\cinf$ function with compact support in $W$, which is
identically equal to $1$ in a neighborhood of $x$, for its pull-back
$\tilde \varphi$ we have:
$$
(2\pi i)^n \ \res \sb 0 (\tilde\omega)= \lim \sb {\delta \to 0} \int
\sb {\{|\tilde f \sb i (z)| = \epsilon \sb i (\delta),\,1
\le i \le n\}} \tilde\varphi \,\cdot\, \tilde\omega = R \sb {\tilde D}
[\tilde\omega] (\tilde\varphi),
$$
from which it follows that
$$
(2\pi i)^n \ \res \sb x (\omega) = R \sb { D}
[\omega] (\varphi).\leqno {(4.9)}
$$
\medskip

\noindent {\bf Remark 4.10.} Given reduced Weil divisors $D_1,\ldots,
D_n$ with finite intersection in a compact $V$-manifold $X$, and a
meromorphic $n$-form $\omega$ whose polar set is contained in the
divisor $D_1\cup\cdots\cup D_n$, it follows from (4.9) and (4.7) that:
$$
(2\pi i)^n \,\sum_{x\in D_1\cap\cdots\cap D_n} \res_x(\omega) =
R_D[\omega](1_X) = (\dbar RP_D[\omega])(1_X) =
-RP_D[\omega](\dbar 1_X) = 0.
$$
This is essentially the proof in [CH, p.~48] of the theorem on the
vanishing of the sum of Grothendieck residues due to Griffiths [G].
\bigskip

{\bf Global Residues.}  We will now generalize the notion of toric
residue to a {\it global residue} defined on an arbitrary
$n$-dimensional compact $V$-manifold $X$.

Given $n+1$ reduced Weil divisors $D_0,\ldots,D_n$ on $X$ such that
$$
D_0\cap\cdots\cap D_n = \emptyset,
$$
the open sets $U_i = X - D_i$ constitute an open cover ${\UU}$ of $X$.
A meromorphic $n$-form $\omega\in\Gamma(X,\widehat\Omega^n_X(*D))$,
with polar set contained in $D=D_0\cup\cdots\cup D_n$, defines a
\v{C}ech cocycle in ${\CC}^n({\UU},\widehat\Omega^n_X)$.  After
passing to the direct limit we obtain a cohomology class $[\omega]\in
H^n(X,\widehat\Omega^n_X)$.  The Dolbeault isomorphism $\eta$ from
Proposition 4.5 (ii) assigns to $[\omega]$ a Dolbeault cohomology
class $\nw\in H^{n,n}_{\dbar}(X)$.

\proclaim Definition 4.11. The {\it global residue} of $\omega$ relative
to the divisors $D_0,\ldots,D_n$ is given by:
$$
\res(\omega) = \biggl(\frac{-1}{2\pi i}\biggr)^n\,\int_X \nw.
$$

For a simplicial toric variety the global residue agrees with the
toric residue.  Indeed, we have already noted in (4.1) that for a
$V$-manifold, our notion of holomorphic forms agrees with the Zariski
differentials and, as shown in [C2, Proposition A.1]:
$${\rm Tr}_X([\omega])=\biggl(\frac{-1}{2\pi i}\biggr)^n\,\int_X \nw.$$

Our next goal is to show that under very mild hypotheses, we can write
the global residue as a sum of local residues.  As above, let
$D_0,\ldots,D_n$ be $n+1$ reduce Weil divisors in $X$ with empty
intersection, and assume that for some $k =0,\ldots,n$, the $n$-fold
intersection
$$
D_{\hat k} = D_0\cap\cdots\cap\widehat{D_k}\cap\cdots\cap D_n
$$
is finite.  If $\omega\in\Gamma(X,\widehat\Omega^n_X(*D))$ and $x\in
D_{\hat k}$ we can write, in a neighborhood of $x$,
$$
\omega = {{\omega'}\over{f_0\cdots f_n}}
$$
where, locally, $\omega'$ is holomorphic and $D_i$ is the support of
$\{f_i = 0\}$, $f_i$ holomorphic.  We will denote by
$\res_{k,x}(\omega)$ the local Grothendieck residue:
$$
\res_{k,x}(\omega) = \res_x\Bigl(
{{\omega'/f_k}\over{f_0\cdots\widehat{f_k}\cdots f_n}}\Bigr).
$$
Note that $x \in D_{\hat k}$ implies that $f_k(x) \ne 0$.

\proclaim Theorem 4.12. If $D_0,\dots,D_n$ are reduced Weil divisors
with empty intersection on a $n$-dimensional compact $V$-manifold $X$,
then for any $\omega\in\Gamma(X,\widehat\Omega^n_X(*D))$, we have:
$$
\res(\omega) =
(-1)^{k}\,  \sum_{x\in D_{\hat k} }\res_{k,x}(\omega)
$$
whenever the intersection $D_{\hat k} =
D_0\cap\cdots\cap\widehat{D_k}\cap\cdots\cap D_n$ is finite.

\noi{\bf Proof.}\ \ There is no loss of generality in assuming $k=n$;
the sign dependence is a consequence of the fact that the global residue
is alternating on the order of the divisors.

The global residue $\res(\omega) = (-1/(2\pi i))^n \int_X \nw$ uses
the Dolbeault isomorphism $\eta$.  However, by Proposition 4.5 (ii),
we can also use the Dolbeault isomorphism ${}'\!\eta$ for currents.
Thus $\res(\omega)$ equals $(-1/(2\pi i))^n$ times the value on the
constant function $1_X$ of any current representing the image under
${}'\!\eta$ of the \v {C}ech cohomology class $[\omega]$.  Hence, the
theorem will follow from the following two assertions:
\medskip
\item{(i)} ${}'\!\eta([\omega])$ is the class of the current
$RP_D[\omega]=R_{D_0,\ldots,D_{n-1}}P_{D_n}[\omega]$.
\item{(ii)} $RP_D[\omega] (1_X) = (2 \pi i)^n
\sum_{x\in D_{\hat n} }\res_{n,x}(\omega)$.
\medskip

Because of the definition of the Dolbeault isomorphism ${}'\!\eta$,
to prove (i), it suffices to construct, for each $i=0,\ldots,n-1$, a \v
{C}ech cochain $\xi^{(i)} \in {\cal C}^i(\UU, \curr n{n-i-1}{X})$
satisfying:
\medskip
\item{(a)} $\delta \xi^{(n-1)} = I(\omega)$ ($\delta$ is the
\v Cech coboundary).
\item{(b)} $\dbar \xi^{(i)} = \delta \xi^{(i-1)}$ for all
$i=1,\ldots,n-1$.
\item{(c)} $\dbar \xi^{(0)} =  R_{D_0,\ldots,D_{n-1}}P_{D_n}[\omega]$.
\medskip
\noi We define
$$
\xi_J^{(n-1)}= \cases{
RP_{D_0}[\omega] & {\rm if} \
$J= \{1, \ldots, n\}$ \cr
0 & {\rm otherwise,} \cr }
$$
and, for any $i= 0, \ldots, n-2$ and any $J \subset \{0, \ldots,n\}$
with cardinality $i+1$,
$$
\xi_J^{(i)}= \cases{
R_{D_0,\ldots,D_{n-i-2}}P_{D_{n-i-1}}[\omega] & {\rm if} \
$J= \{n-i, \ldots,n\}$ \cr
 0    & {\rm otherwise.} \cr }
$$
It is understood that the above currents $RP_{D_0}[\omega]$ and
$R_{D_0,\ldots,D_{n-i-2}}P_{D_{n-i-1}}[\omega]$ are restricted to the
appropriate open sets $U_J=\dsize\cap_{j\in J} U_j$.  We will
generally not indicate the restriction when it is irrelevant or clear
from the context.

To verify (a), note that $\delta \xi^{(n-1)}$ is the cochain assigning
to $U_0\cap\cdots\cap U_n$ the current $RP_{D_0}[\omega]$ restricted
to this open set.  Since $D_0$ is disjoint from $U_0\cap\cdots\cap
U_n$, the definition of $RP_{D_0}[\omega]$ implies that it must agree
with $I(\omega)$.

Suppose now that $1\leq i \leq n-1$, then it follows from (4.7) that
$$
(\dbar \xi^{(i)})_J =  \cases{
R_{D_0,\ldots,D_{n-i-1}}[\omega] & {\rm if} \
$J= \{n-i, \ldots, n\}$ \cr
 0    & {\rm otherwise.} \cr }
$$
On the other hand, clearly $(\delta \xi^{(i-1)})_J =0$ if
$J$ is not an index set of the form $J=J_j = \{j,n-i+1,\ldots,n\}$ for
some $j=1,\ldots,n-i$.  But, if $j<n-i$, then $(\delta
\xi^{(i-1)})_{J_j}$ also vanishes---as a consequence of
(4.6)---since it is the restriction to the open set $U_{J_j}\subset
U_j$ of the current $R_{D_0,\ldots,D_{n-i-1}}P_{D_{n-i}}[\omega]$ and
$j<n-i$.

It remains to consider the case $J=\{n-i,\ldots,n\}$.  Then, $(\delta
\xi^{(i-1)})_{J}$ is the restriction to $U_J$ of
$R_{D_0,\ldots,D_{n-i-1}}P_{D_{n-i}}[\omega]$.  But, since $U_J\subset
U_{n-i}$, we deduce, again from (4.6), that
$$
R_{D_0,\ldots,D_{n-i-1}}P_{D_{n-i}}[\omega]
\lower.04truein\hbox{$\big\vert$}
\lower.07truein\hbox{$\scriptstyle{U_J}$} =
R_{D_0,\ldots,D_{n-i-1}}[\omega]\lower.04truein\hbox{$\big\vert$}
\lower.07truein\hbox{$\scriptstyle{U_J}$}.
$$
Thus, (b) is satisfied.

The final assertion (c) is proved in a similar way: the cochain $\dbar
\xi^0$ assigns the zero current to the open sets $U_j$, $j<n$ and the
residue current $R_{D_0,\ldots,D_{n-1}}[\omega]$ to $U_n$.  But, then,
it follows from (4.6) that $\dbar \xi^0$ agrees with the global
current $RP_D[\omega]$.

The verification of (ii) now reduces to the local formula (4.9).
Indeed, since the support of the principal value $RP_D[\omega]$ is
contained in the finite set $D_{\hat n} = D_0\cap\cdots\cap D_{n-1}$,
its value on the constant function $1_X$ is the same as the value on
any function $\psi$ which is equal to one on a neighborhood of each of
the points in $D_{\hat n}$.  We may choose such a function $\psi$ of
the form $\psi = \sum_{x\in D_{\hat n}} \psi_x$, where $\psi_x$ is
equal to $1$ in a neighborhood of $x$ and the supports of the
$\psi_x$'s are mutually disjoint and disjoint from $D_n$ as well. Then
$$
RP_D[\omega](1_X) = \sum_{x\in D_{\hat
n}}RP_D[\omega]\lower.04truein\hbox{$\big\vert$}
\lower.07truein\hbox{$\scriptstyle{U_n}$}(\psi_x) =
\sum_{x\in D_{\hat n}}R_{D_0,\ldots,
D_{n-1}}[\omega]\lower.04truein\hbox{$\big\vert$}
\lower.07truein\hbox{$\scriptstyle{U_n}$}(\psi_x)
$$
where the last equality follows from (4.6).  But now, (4.9) yields
$$
RP_D[\omega](1_X) =(2\pi i)^n\sum_{x\in D_{\hat n}}\res_x(\omega).
\ \ \diamond
$$
\smallskip

\noi {\bf Remarks 4.13.}\ (i)\ To understand why we need currents in
the proof of Theorem~4.12, we will sketch a proof for the case $n=2$
using forms rather than currents.  The argument will be less than
rigorous.

We have $D_0\cap D_1\cap D_2 = \emptyset$ in $X$.  Let $T_j(\epsilon)$
be a fundamental system of (open) tubular neighborhoods of $D_j$, and
let $S_j(\epsilon) = \partial T_j(\epsilon)$ and $E_j(\epsilon) = X-
{{T_j(\epsilon)}}$.  Also, for $i,j,k$ distinct indices from $0$ to
$2$, consider the intersections $C_{ijk}(\epsilon) = E_i(\epsilon)
\cap E_j(\epsilon) \cap E_k(\epsilon)$, $C_{ij}(\epsilon) =
E_i(\epsilon) \cap E_j(\epsilon) \cap S_k(\epsilon)$ and
$C_{i}(\epsilon) = E_i(\epsilon) \cap S_j(\epsilon) \cap
S_k(\epsilon)$.  We will assume that these sets are homology chains of
(real) codimension $0,1$, and $2$ respectively and that their
boundaries behave as one would expect.

Next recall the procedure to define $\nw$.  Let $\{\sigma_0, \sigma_1,
\sigma_2\}$ be a partition of unity subordinated to the covering
$\UU$.  Then, beginning with $\omega\in \Gamma(U_0\cap U_1\cap U_2,
\widehat\Omega^2_X)$, define $\xi_{ij} = (-1)^k \,\sigma_k\,\omega \in
\Gamma(U_i\cap U_j, {\cal E}^{2,0}_X)$, which implies $\omega = \delta
(\xi_{ij}) = \xi_{12} - \xi_{02} + \xi_{01}$.  Next, define $\xi_i =
\pm \sigma_j \,\dbar \xi_{ij} \pm \sigma_k\,\dbar \xi_{ik}$, with the
signs chosen so that $\dbar \xi_{ij} = \delta(\xi_i) = \xi_j - \xi_i$.
Finally, $\nw$ is defined to be the global $(2,2)$-form $\dbar
\xi_i$ in $U_i$.

To compute the global residue $(2\pi i)^{-2}\int_X\nw$, we first observe
$$
\int_X\,\nw = \lim_{\epsilon\to 0}\,\int_{C_{012}(\epsilon)} \,\nw =
\lim_{\epsilon\to 0}\,\int_{C_{012}(\epsilon)} \, \dbar \xi_0.
$$
Since $\xi_0$ has bidegree $(2,1)$, $d\xi_0 = \dbar \xi_0$, so we can
apply Stokes' Theorem to write
$$
\int_{C_{012}(\epsilon)} \, \dbar \xi_0  =
\int_{C_{12}(\epsilon)} \,  \xi_0\, +
\int_{C_{02}(\epsilon)} \,  \xi_0\, +
\int_{C_{01}(\epsilon)} \,  \xi_0.
$$
On the other hand, $\xi_0 = \pm\, \sigma_1\,\dbar \sigma_2\wedge\omega
\pm\, \sigma_2\,\dbar \sigma_1\wedge\omega$, and therefore $\xi_0 = 0$
in $S_1(\epsilon)$ and $S_2(\epsilon)$ for sufficiently small
$\epsilon$.  Consequently,
$$
\int_X\,\nw =
 \lim_{\epsilon\to 0}\,\int_{C_{12}(\epsilon)} \,  \xi_0 =
 \lim_{\epsilon\to 0}\,\int_{C_{12}(\epsilon)} \, \xi_1 -
\dbar\xi_{01}.
 $$
Once again, $\int_{C_{12}(\epsilon)} \, \xi_1
$vanishes for $\epsilon$ sufficiently small and, using Stokes' Theorem,
we write
$$
\int_X\,\nw = \lim_{\epsilon\to 0}\,\int_{C_{2}(\epsilon)} \,
\xi_{01} + \lim_{\epsilon\to 0}\,\int_{C_{1}(\epsilon)} \,
\xi_{01}.
$$
Since $\xi_{01} = \sigma_2\,\omega$ vanishes on $S_2(\epsilon)$
for $\epsilon$ sufficiently small, we have
$$
\int_X\,\nw = \lim_{\epsilon\to 0}\,\int_{C_{2}(\epsilon)} \,
\sigma_2\,\omega.
$$
Finally, for $\epsilon$ sufficiently small, $\sigma_2$ is identically
$1$ in $S_0(\epsilon)\cap S_1(\epsilon)$, so that
$$
\int_X\,\nw = \lim_{\epsilon\to 0}\,\int_{C_{2}(\epsilon)} \,
\omega = (2\pi i)^2  \sum_{x\in D_0\cap D_1 }\res_{x}(\omega),
$$
which gives the desired formula for the global residue.

The use of residual currents in making the above argument rigorous is
twofold: first of all, the local nature of the residual currents
definition obviates the need to construct global cycles of
integration---a step which is not always possible; moreover, the
concept of admissible paths explains the passage to the limit
necessary for the vanishing of the various integrals.

\medskip
\noi {(ii)}\ For an example of how
$D_0,\ldots,\widehat{D_k},\ldots,D_n$ can fail to satisfy the
finiteness condition in Theorem 4.12, let $X =
\P^1\times\P^1$, and consider the divisors $D_0 =
\{0\}\times\P^1$, $D_1 = (\{\infty\}\times\P^1) \cup
(\P^1\times\{\infty\})$ and $D_2 = (\{\infty\}\times\P^1) \cup
(\P^1\times\{0\})$.  Then $D_0\cap D_1\cap D_2 = \emptyset$, yet
$D_1\cap D_2$ is infinite since it contains $\{\infty\}\times\P^1$.
\medskip

	In light of the last remark, it would be useful to know when
the $n$-fold intersection $D_{\hat k} =
D_0\cap\cdots\cap\widehat{D_k}\cap\cdots\cap D_n$ is finite.  Here
is one criterion.

\proclaim Lemma 4.14. Let $D_0,\ldots,D_n$ be reduced Weil divisors
with empty intersection on a $n$-dimensional projective variety $X$.
If $D_k$ is the support of an ample divisor, then $D_{\hat k}$ is
finite.

\noi {\bf Proof.} This is immediate since $D_{\hat k}\cap D_k =
\emptyset$ implies that $D_{\hat k}$ is a complete subvariety of the
affine variety $X - D_k$.$\ \ \diamond$
\medskip

	When applied to toric residues, these results yield Theorem
0.4 which, in turn, may be used to give new proofs of some basic
results concerning toric residues.  We will conclude this section with
three such applications of Theorem 0.4:
\medskip

{\bf Res = $\pm$1.}\ The first application is an alternate proof of
Proposition 2.4 when the toric variety $X$ is simplicial.  We resume
the notation of \S 2, where we have a $n$-dimensional cone $\sigma$
and the variables are labelled $x_1,\ldots,x_n$ (corresponding to the
generators of $\sigma$) and $z_1,\ldots,z_r$ (corresponding to the
other generators).  Our goal is to prove that
$$
\res\biggl(\frac{\Omega}{(z_1\cdots z_r)\cdot x_1\cdots x_n}\biggr) =
\pm1
$$
using Theorem 0.4.  Since the divisors
$$
D_0 = \{z_1\cdots z_r =0\},\quad D_i=\{x_i = 0\},\ i=1,\ldots,n
$$
have empty intersection and $D_1\cap\cdots\cap D_n = \{p\}$ is a
single point, the hypotheses of Theorem~0.4 are satisfied, so that
$$
\res\biggl(\frac{\Omega}{(z_1\cdots z_r)\cdot x_1\cdots x_n}\biggr) =
\res_p(\omega_{\sigma}),
$$
where $\omega_{\sigma}$ is the restriction of $\Omega/\bigl(
(z_1\cdots z_r)\cdot x_1\cdots x_n\bigr)$ to the affine open set
$X_{\sigma}\subset X$.  But we have seen that $X_{\sigma}\simeq
\C^n/G(\sigma)$, and since $G(\sigma)$ is a small subgroup, this
defines a standard model.  Moreover, as noted in (2.9), the pullback
to $\C^n$ of $\omega_{\sigma}$ is given by
$$
\tilde\omega_{\sigma} = \pm{|G(\sigma)|\, dx_1\wedge\cdots\wedge dx_n
\over x_1\cdots x_n}.
$$
Therefore
$$
\res_p(\omega_{\sigma}) = \frac 1{|G(\sigma)|}
\res_0(\tilde\omega_{\sigma}) = \pm\,\res_0\biggl(
{dx_1\wedge\cdots\wedge dx_n \over x_1\cdots x_n}\biggr)=\pm 1,
$$
which proves Proposition 2.4 when $X$ is simplicial.
\bigskip

{\bf Sums of Residues in a Torus.}\ Let $f_1,\ldots,f_n$ be
$n$-variate Laurent polynomials with a finite set of common zeroes $Z
= Z(f_1,\ldots,f_n)$ in the torus $T=(\C^*)^n$.  Given a Laurent
polynomial $q$, we get the differential form
$$
\phi = \frac q{f_1\cdots f_n}\,\frac {dt_1}{t_1}\wedge
\cdots\wedge \frac {dt_n}{t_n}.
$$
The operator which assigns to $q$ the sum of local residues
$\sum_{x\in Z} \res_x(\phi)$ has interesting applications in a number
of different contexts.  In certain cases, it is possible to use
Theorem~0.4 to give a global interpretation of this sum.

We assume that there exists a simplicial toric compactification $X$ of
$T$ such that if $D_i$ is the closure in $X$ of the hypersurface
$\{f_i = 0\} \subset (\C^*)^n$ and $D_0 = X - T$ is the ``divisor at
infinity'', then
$$
D_0\cap D_1\cap \cdots \cap D_n = \emptyset.
$$
Such a (smooth) compactification exists, for example, if the
polynomials $f_i$ are nondegenerate in the sense of Khovanskii [K1].

In this situation, the meromorphic form $\phi$ has an extension to $X$
which can be written as
$$
\Phi  = \frac {Q\,\Omega}{F_0\cdots F_n},
$$
where $Q,F_0,\ldots,F_n$ are homogeneous polynomials in the coordinate
ring of $X$ such that $D_i = \{F_i=0\}$, and $\Omega$ is the Euler
form of $X$.  Then it follows from Theorem 0.4 that
$$
\sum_{x\in Z} \res_x(\phi) = \res_F(Q).
$$

If we assume, in addition, that the Newton polyhedron of $q$ is
contained in the interior of the Minkowski sum of the Newton polyhedra
corresponding to $f_1,\ldots,f_n$, then one may show that $Q$ is a
multiple of $F_0$ and hence $\res_F(Q)=0$ which gives the classical
Euler-Jacobi Theorem in this setting [K2].  In fact, as in this case
$\Phi$ has poles only on the union of the $n$ divisors $D_1 \cup
\cdots \cup D_n$, whose intersection is contained in the torus, the
vanishing of the sum of the local residues of $\phi$ follows directly
from the result of Griffiths recalled in Remark 4.10.
\bigskip

{\bf Toric Jacobians.}\ For our third application, we use Theorem~0.4
to give an alternate proof of Theorem 5.1 (ii) of [C2] for a
simplicial toric variety.  This result asserts that the {\it toric
Jacobian} $J$ of $F_0,\ldots,F_n \in S_\alpha$ (as defined in [C2,
Proposition 4.1]) has nonzero toric residue.  More precisely, if
$\alpha$ is ample and the $F_i$ don't vanish simultaneously on $X$,
then we will show that the equality
$$
\res_F(J) = (D^n)
$$
follows from Theorem 0.4.  Here, $(D^n)$ is the $n$-fold
intersection number of any divisor $D$ with $[D] = \alpha$.  Note
that $J \in S_{\rho}$, where $\rho = (n+1)\alpha - \beta$ is the
critical degree for the $F_i$.

To prove this, let $\omega_F(J) = J\,\Omega/(F_0\cdots F_n)$.  Then
Theorem 0.4 implies
$$
\res_F(J) = \sum_{x\in D_{\hat 0}
}\res_{0,x}(\omega_F(J)) = \sum_{x\in D_1\cap\cdots\cap D_n}
\res_x\biggl({(J/F_0)\,\Omega \over F_1\cdots F_n}\biggr).
\leqno{(4.15)}
$$
We will show that each local residue $\res_{0,x}(\omega_F(J))$ is a
local intersection multiplicity of $D_1,\ldots,D_n$ at $x$, which will
prove that $\res_F(J)$ is the intersection number $(D_1\cdots D_n) =
(D^n)$.

Given $x \in D_{\hat 0}$, let $\sigma$ be a $n$-dimensional cone such
that $x$ lies in the affine open set $X_{\sigma}$.  Since $\sigma$ is
simplicial, we write the variables as $x_1,\ldots,x_n,z_1,\ldots,z_r$.
Then the form $\omega_F(J)$, restricted to $X_\sigma$, may be written
in appropriate coordinates as:
$$
\frac{k(x_1,\ldots,x_n)\, dx_1\wedge\cdots\wedge dx_n}{f_0\cdots f_n}
$$
where
$$
k(x_{1},\ldots,x_{n})={\rm det}\
\pmatrix{
 f_0&\cdots& f_n\cr
\partial  f_0/\partial x_{1}&\cdots&\partial
 f_n/\partial x_{1}\cr \vdots & & \vdots \cr
\partial  f_0/\partial x_{n}&\cdots&\partial
 f_n/\partial x_{n}\cr}
$$
and $f_i(x_{1},\ldots,x_{n})$ is the function obtained from
$F_i(x_1,\ldots,x_n,z_1,\ldots,z_r)$ by setting $z_j =1$ for
$j=1,\ldots,r$.  It follows that
$$
\res_{0,x}(\omega_F(J)) = \res_x\biggl(\frac
{(k/f_0)\,dx_1\wedge\cdots\wedge dx_n}{f_1\cdots f_n}\biggr).
$$
However, expanding the determinant for $k$ along the first row and
using $f_0(x) \ne 0$, we see that
$$
k/f_0 \equiv \det(\partial f_i/\partial x_j: 1\le i,j\le n) \bmod
\langle f_1,\ldots,f_n\rangle
$$
in the local ring $\OO_{X,x}$.  Consequently,
$$
\res_{0,x}(\omega_F(J)) = \res_x\biggl(\frac
{\det(\partial f_i/\partial
x_j)\,dx_1\wedge\cdots\wedge dx_n}{f_1\cdots f_n}\biggr),
$$
and this last residue equals the local intersection multiplicity of
$D_1,\ldots,D_n$ at $x$ (this is well-known in the smooth case and is
easy to prove for $V$-manifolds).  By (4.15), it follows that the
toric residue of the toric jacobian equals the intersection number
$(D_1\cdots D_n)$.
\bigskip

\noi {\bf Remarks 4.16.}\ (i)\ Since all the divisors $D_i$ have the
same degree, we can write $\res_F(J) = (D_0\cdots \widehat{D_k} \cdots
D_n)$ for any $k=0,\ldots,n$.
\medskip
\noi (ii)\ The intersection number $(D_1\cdots D_n)$ can also be
interpreted as the degree of the map $F = (F_0,\ldots,F_n) : X \to
\P^n$ (see [C2, Theorem 5.1] for a careful proof).  Thus the toric
Jacobian has the property that its toric residue is given by
$\res_F(J) = \deg(F)$.  This will be useful in \S5.
\bigskip
\bigskip

\noi {\bf \S5. Toric Residues as Point Residues in the Equal Degree
Case}
\bigskip

When $X = \P^n$ and $F_0,\ldots,F_n$ all have the same degree, the
toric residue $\res_F(H)$ equals the classical Grothendieck residue at
$0 \in \C^{n+1}$, i.e.,
$$
\res_F(H) = {1\over (2\pi i)^{n+1}}\int_{|F_i| = \epsilon}
{H\,dx_0\wedge\cdots \wedge dx_n \over F_0\cdots F_n}
$$
(see [PS, 12.10]).  Thus, in the projective case, the toric residue is
a point residue computed on the related space $\C^{n+1}$.  In Theorem
5.8 below, we will generalize this result to a complete simplicial
toric variety $X$, assuming that the $F_i$ have the {\it same} degree
$\alpha$ in the homogeneous coordinate ring $S$.

We first describe the space we will use for computing toric residues
on $X$.  Given $\alpha \in A_{n-1}(X)$, let $S_{*\alpha} =
\oplus_{k\ge0} S_{k\alpha}$ and set
$$
X_\alpha = {\rm Spec}(S_{*\alpha}).
$$
Note that the natural grading of $S_{*\alpha}$ induces a $\C^*$ action
on $X_\alpha$.

\proclaim Proposition 5.1. If $X$ is a complete simplicial toric variety
and $\alpha$ is ample, then $X_\alpha$ has the natural structure of an
affine toric variety.  Furthermore, if $0 \in X_\alpha$ is the unique
fixed point of the torus action, then $X_\alpha - \{0\}$ is simplicial
and $\C^*$ acts on $X_\alpha - \{0\}$ with finite stabilizers and $X$
as geometric quotient.

\noi{\bf Proof.} Consider $\R\oplus N_\R$ with the lattice $\Z\oplus
N$.  Elements of $\R\oplus N_\R$ will be written $\lambda e_0 + v$,
where $\lambda \in \R$ and $v \in N_\R$.  Now let $D = \sum_i a_i D_i$
(where $\sum_i$ denotes $\sum_{i=1}^{n+r}$) be a divisor on $X$ whose
class is $\alpha$, and let $\psi : N_\R \to \R$ be its support
function.  This means $\psi(\eta_i) = -a_i$, where the $\eta_i$
generate the 1-dimensional cones of the fan of $X$.  Given this data,
let $\widetilde\sigma \subset \R\oplus N_\R$ be the cone generated by
the vectors $\tilde \eta_i = a_i e_0 + \eta_i$.  Equivalently,
$\widetilde\sigma$ is generated by the graph of $-\psi$ in $\R\oplus
N_\R$.  Since $\psi$ is strictly upper convex ($D$ is ample), we see
that $\widetilde\sigma$ is a strongly convex rational polyhedral cone.

We next observe that the semigroup ring $\C[\widetilde\sigma^\vee \cap
(\Z\oplus M)]$ is naturally isomorphic to $S_{*\alpha}$.  To prove
this, first note that
$$
\eqalign{ke_0+m \in \widetilde\sigma^\vee\cap(\Z\oplus M) &\iff
	\langle ke_0+m,\tilde \eta_i\rangle \ge 0\ \hbox{for all $i$}\cr
&\iff \langle m,\eta_i\rangle + k\,a_i \ge 0\ \hbox{for all $i$}\cr
&\iff \Pi_i x_i^{\langle m,\eta_i\rangle + k\,a_i} \in
S_{k\alpha}\cr}
$$
(where $\Pi_i$ denotes $\Pi_{i=1}^{n+r}$).  Since all monomials in
$S_{k\alpha}$ can be described in this form (see \S1 of [C1]), the
observation follows easily.  Thus $X_\alpha$ is an affine toric
variety.

The torus action on $X_\alpha$ has a unique fixed point which we
denote by $0$.  Furthermore, the complement $X_\alpha - \{0\}$ is the
toric variety whose fan is the boundary of $\widetilde\sigma$.  This
fan is the graph of $-\psi$, so the strict convexity of $\psi$ implies
that under the projection $\pi: \R\oplus N_\R \to N_\R$, each cone of
the boundary fan maps naturally to the corresponding cone in the fan
of $X$.  Thus the projection $\pi$ induces a map of toric varieties
$X_\alpha - \{0\} \to X$.  We leave to the reader the straightforward
proof that $X_\alpha - \{0\}$ is simplicial since $X$ is.

Since $X$ is simplicial, we can write
$$
X = (\C^{n+r} - Z)/G, \leqno{(5.2)}
$$
where $G = {\rm Hom}_\Z(A_{n-1}(X),\C^*)$ and the exceptional set $Z$
is a union of coordinate subspaces determined by the fan of $X$ (see
[BC, Theorem 1.9]).  The correspondence $\eta_i \leftrightarrow \tilde
\eta_i$ implies that $X$, $X_\alpha$ and $X_\alpha - \{0\}$ have the
same homogeneous coordinate ring (though the gradings may differ), and
the map $X_\alpha - \{0\} \to X$ shows that the fans of $X$ and
$X_\alpha - \{0\}$ are combinatorially equivalent.  Thus $X$ and
$X_\alpha - \{0\}$ have the same exceptional set $Z$.  Hence
$$
X_\alpha - \{0\} = (\C^{n+r} - Z)/H, \leqno{(5.3)}
$$
where $H = {\rm Hom}_\Z(A_n(X_\alpha - \{0\}),\C^*)$.  To compare $G$
and $H$, we use the commutative diagram
$$
\matrix{0 & \to & M & \to & \Z^{n+r} & \to & A_{n-1}(X) & \to & 0\cr
&& \downarrow && \parallel && \downarrow && \cr
0 & \to & \Z\oplus M & \to & \Z^{n+r} & \to & A_{n}(X_\alpha - \{0\})
& \to & 0\cr}
$$
to conclude that we have an exact sequence
$$
0 \longrightarrow \Z \longrightarrow A_{n-1}(X) \longrightarrow
A_n(X_\alpha - \{0\}) \longrightarrow 0,
$$
where $1 \in \Z$ maps to $\alpha \in A_{n-1}(X)$.  Applying ${\rm
Hom}_\Z(-,\C^*)$, we can identify $H$ with the subgroup $\{g \in G :
g(\alpha) = 1\} \subset G$, so that $g \mapsto g(\alpha)$ induces an
isomorphism $G/H \simeq \C^*$.

	Comparing (5.2) and (5.3), $X$ is the quotient of $X_\alpha
-\{0\}$ by $G/H \simeq \C^*$.  Furthermore, the proof of Theorem 1.9
of [BC] shows that the $G$-action in (5.2) has finite stabilizers, and
it follows that the $\C^*$-action on $X_\alpha$ must also have finite
stabilizers.  To describe this action more explicitly, note that $G$
acts on $S_{k\alpha}$ by $g\cdot F = g(k\alpha)\, F = g(\alpha)^k F$.
Since $H$ acts trivially by definition, the action of $G/H \simeq
\C^*$ is exactly the action that gives the grading of $S_{*\alpha}$.
This completes the proof of the proposition.$\ \
\diamond$
\bigskip

\noi {\bf Remarks 5.4.}\ (i)\ When $\alpha$ is very ample (always true
when $X$ is smooth), then $X_\alpha$ is the affine cone of $X$ in the
projective embedding given by $\alpha$.
\medskip
\noi (ii)\ Besides being a geometric quotient, the map $X_\alpha -
\{0\} \to X$ is a {\it combinatorial quotient\/} in the sense of [KSZ,
p.~645].
\medskip
\noi (iii)\ If we add the 1-dimensional cone generated by $e_0$ to
$\widetilde\sigma$ and subdivide accordingly, we get a toric variety
$\widetilde X_\alpha$ which maps naturally to $X$.  In [R, Section 3],
it is proved that $\widetilde X_\alpha \to X$ is the total space of
the line bundle ${\cal O}_X(-\alpha)$.  Thus $\widetilde X_\alpha \to
X_\alpha$ is a blow-up of $0 \in X_\alpha$ with exceptional fiber
isomorphic to $X$.  Conversely, we can view $X_\alpha$ as the variety
obtained by blowing down the zero section of ${\cal O}_X(-\alpha)$.
\medskip
\noi (iv)\ Although $X_\alpha-\{0\}$ is simplicial, $0 \in X_\alpha$
can be very singular.  For example, let $X = \P^1\times\P^1$ and
$\alpha = (1,1)$.  The coordinate ring for $X$ is $S = \C[x,y,z,t]$,
where $\deg(x) = \deg(y) = (1,0)$ and $\deg(z) = \deg(t) = (0,1)$.
Then $X_\alpha$ is the singular affine hypersurface defined by $A\,D -
B\,C = 0$ in $\C^4$ since this hypersurface is the affine cone over
the Segre embedding $\P^1\times\P^1 \hookrightarrow \P^3$.  Note that
$X_\alpha$ is not simplicial at the origin.
\bigskip

We next discuss differential forms on $X$ and $X_\alpha$.  As we saw
in the proof of Proposition~5.1, $X$ and $X_\alpha$ have the same
homogenous coordinate ring (though graded differently).  By (2.8), $X$
has the Euler form
$$
\Omega = \sum_{|I| = n} \det(\eta_I)\,\hat x_I\, dx_I.
$$
Now let $\sum_i a_i\,D_i$ be a divisor in the class of
$\alpha$ and consider the $(n+1)$-form
$$
\Omega_\alpha = \Big(\sum_i a_i {dx_i \over
x_i}\Big) \wedge \Omega.
$$

\proclaim Lemma 5.5. Let $\beta = \sum_i \deg(x_i)\in A_{n-1}(X)$ and
$\rho = (n+1)\alpha - \beta$.  Then:
\vskip0pt
\item{(i)} $\Omega_\alpha$ is the Euler form of $X_\alpha$.
\item{(ii)} If $\theta$ is any Euler vector field for $X$ (which can
be regarded as a map $\theta : A_{n-1}(X) \to \C$), we have
$$
\theta\ip\Omega_\alpha = \theta(\alpha)\,\Omega.
$$
\item{(iii)} If $J \in S_\rho$ is the toric Jacobian of $F_0,\dots,F_n
\in S_\alpha$ (see [C2, \S4]), then
$$
J\,\Omega_\alpha = dF_0\wedge\cdots\wedge dF_n.
$$
\vskip0pt

\noi{\bf Proof.} To define $\Omega$, we used a basis $m_1,\dots,m_n$
of $M$.  Then $e_0$ and $m_j$ for $j > 0$ form a basis of $\Z\oplus
M$, and from the proof of Proposition 9.5 of [BC] (which
is easily seen to hold in the non-simplicial case), we see that the
Euler form of $X_\alpha$ is
$$
x_1\cdots x_{n+r} \Big({dt_0\over t_0}\wedge\cdots\wedge{dt_n\over
t_n}\Big), \leqno{(5.6)}
$$
where $t_0 = \Pi_i x_i^{\langle e_0,\tilde \eta_i\rangle} = \Pi_i
x_i^{a_i}$ and $t_j = \Pi_i x_i^{\langle m_j,\tilde \eta_i\rangle} =
\Pi_i x_i^{\langle m_j,\eta_i\rangle}$ for $j > 0$.  Since $dt_0/t_0 =
\sum_i a_i dx_i/x_i$ and $\Omega = x_1\cdots x_{n+r}
(dt_1/t_1\wedge\cdots\wedge dt_n/t_n)$ (also by Proposition 9.5 of
[BC]), we see that $\Omega_\alpha$ is the Euler form of $X_\alpha$.

For the second part of the lemma, first note that $\theta\ip\Omega =
0$ by Lemma 6.2 of [C2].  Thus
$$
\theta\ip\Omega_\alpha =
\theta\ip\Big(\big(\textstyle{\sum_i} a_i
{dx_i \over x_i}\big) \wedge \Omega\Big)
= \Big(\theta\ip\textstyle{\sum_i} a_i
{dx_i \over x_i}\Big) \cdot \Omega.
$$
However, if $\theta = \sum_i b_i\, x_i\, \partial/\partial x_i$, then
$\theta\ip\sum_i a_i dx_i/x_i = \sum_i a_i b_i = \theta(\alpha)$,
which gives the desired formula.  (For more background on Euler vector
fields, see 3.8--3.10 of [BC].)

Turning to the final part of the lemma, note that each $F_i$ lies in
$S_\alpha$ and hence gives a function on $X_\alpha = {\rm
Spec}(S_{*\alpha})$.  Further, the functions $t_0,\dots,t_n$ introduced
above are coordinates on the torus $T_{X_\alpha} \subset X_\alpha$.
Thus, if we restrict $F_i$ to the torus, we can write $F_i =
\widetilde F_i(t_0,\dots,t_n)$.  Then
$$
dF_0\wedge\cdots\wedge dF_n = \det(\partial \widetilde F_i/\partial
t_j)\,dt_0\wedge\cdots\wedge dt_n.
$$
Comparing this to the formula (5.6) for $\Omega_\alpha$, we see that
$$
dF_0\wedge\cdots\wedge dF_n = \widetilde J\,\Omega_\alpha
$$
for some rational function $\widetilde J$.

It remains to show that $\widetilde J$ is the toric Jacobian $J$ from
[C2].  Pick an Euler formula $\theta$ such that $\theta(\alpha) \ne
0$.  We can find such a $\theta$ since $\alpha$ is ample and hence has
infinite order in $A_{n-1}(X)$ (see also Lemma 10.5 of [BC]).  Then,
by (ii) and the above equation for $\widetilde J$,
$$
\eqalign{\theta(\alpha)\cdot \widetilde J\,\Omega =
\theta\ip(\widetilde J\,\Omega_\alpha) &=
\theta\ip(dF_0\wedge\cdots\wedge dF_n)\cr
&= \textstyle{\sum_{i=0}^n} (-1)^i (\theta\ip
dF_i)\,dF_0\wedge\cdots\wedge \widehat{dF_i} \wedge\cdots\wedge dF_n\cr
&= \theta(\alpha)\cdot \textstyle{\sum_{i=0}^n} (-1)^i F_i\,
dF_0\wedge\cdots\wedge \widehat{dF_i} \wedge\cdots\wedge dF_n,\cr}
$$
where the last equality follows because $\theta\ip dF =
\theta(\alpha)\,F$ for all $F \in S_\alpha$.  However, on the bottom
row, the expression on the right equals $\theta(\alpha) \cdot
J\,\Omega$ by [C2].  Then $\widetilde J = J$ follows since
$\theta(\alpha) \ne 0$, and (iii) is proved.$\ \ \diamond$
\medskip

	Given $F_0,\dots,F_n \in S_\alpha$, we next consider the
integral on $X_\alpha$
$$
\int_{\{|F_i| = \epsilon, 0\le i \le n\}} {H\,\Omega_\alpha \over
F_0\cdots F_n}, \leqno{(5.7)}
$$
where $\epsilon > 0$, the cycle $\{|F_i| = \epsilon, 0 \le i \le n\}$
is oriented using $d({\rm arg}\,F_0)\wedge\cdots\wedge d({\rm
arg}\,F_n)$, and $H \in S_\rho$ for $\rho = (n+1)\alpha - \beta$.  To
make sense of (5.7), first note that $H\,\Omega_\alpha/(F_0\cdots F_n)$
is a meromorphic form on the $V$-manifold $X_\alpha - \{0\}$.
Furthermore, each $F_i$ is a polynomial function on $X_\alpha$ and
$\{|F_i| = \epsilon, 0 \le i \le n\} \subset X_\alpha - \{0\}$.  It
follows that (5.7) exists whenever $(F_0,\dots,F_n) : X_\alpha \to
\C^{n+1}$ is finite.  We can now state the main result of this
section.

\proclaim Theorem 5.8.  Assume that $X$ is complete and simplicial,
$\alpha$ is ample, and $F_0,\ldots,F_n \in S_\alpha$ don't vanish
simultaneously on $X$.  Then:
\vskip0pt
\item{(i)} The map $(F_0,\dots,F_n) : X_\alpha \to \C^{n+1}$ is
finite.
\item{(ii)} If $\rho = (n+1)\alpha - \beta$ is the
critical degree of $F_0,\dots,F_n$, then for every $H \in S_\rho$,
$$
\res_F(H) = {1\over(2\pi i)^{n+1}} \int_{\{|F_i| = \epsilon, 0\le
i\le n\}} {H\,\Omega_\alpha \over F_0\cdots F_n}.
$$

\noi{\bf Proof.}  By [C2, Proposition 3.2], we know that
$S_{*\alpha}/\langle F_0,\ldots,F_n\rangle$ has finite dimension over
$\C$, so that by definition, $F_0,\ldots,F_n$ is a homogeneous system
of parameters for $S_{*\alpha}$.  It follows from [BH, Theorem 1.5.17]
that $S_{*\alpha}$ is finitely generated as a module over the subring
$\C[F_0,\ldots,F_n]$.  Thus $\widetilde F = (F_0,\ldots,F_n) : X_\alpha
\to \C^{n+1}$ is finite, which proves (i).

	To prove (ii), we first observe that each side of the identity
in (ii) vanishes when $H \in \langle F_0,\ldots,F_n\rangle$.  This is
obviously true for the toric residue, and for the integral (5.7), one
uses the usual argument (see [GH, pp.~650--651]).  Since we know
$S_\rho/\langle F_0,\ldots,F_n\rangle_\rho$ is one dimensional and the
toric Jacobian $J$ has nonzero toric residue (see \S4), it suffices to
check that (ii) holds for $H = J$.

	By Remark 4.16, we know that $\res_F(J) = \deg(F)$, where $F =
(F_0,\ldots,F_n)$, regarded as a map $F : X \to \P^n$.  On the other
hand, by Lemma 5.5, we have
$$
{J\,\Omega_\alpha\over F_0\cdots F_n} = {dF_0\wedge\cdots\wedge dF_n
\over F_0\cdots F_n} = \widetilde F^*\biggl({dz_0\wedge\cdots\wedge dz_n
\over z_0\cdots z_n}\biggr),
$$
where $z_0,\ldots,z_n$ are coordinates on $\C^{n+1}$ and $\widetilde F
= (F_0,\ldots,F_n)$, now regarded as a map $\widetilde F : X_\alpha
\to \C^{n+1}$.  It follows that
$$
{1\over(2\pi i)^{n+1}} \int_{\{|F_i| = \epsilon, 0\le i\le n\}}
{J\,\Omega_\alpha \over F_0\cdots F_n} = {\deg(\widetilde F)\over(2\pi
i)^{n+1}} \int_{\{|z_i| = \epsilon, 0\le i\le
n\}}\!\!{dz_0\wedge\cdots\wedge dz_n \over z_0\cdots z_n} =
\deg(\widetilde F)
$$
since $\widetilde F$ is finite by (i).

	Thus, to prove (ii) for $J$, we must show that $\deg(F) =
\deg(\widetilde F)$.  However, as noted in the proof of Proposition 5.1,
the $\C^*\simeq G/H$ action on $X_\alpha$ satisfies $g\cdot F_i =
g(\alpha)\,F_i$ for $g \in G$.  It follows that $\widetilde F :
X_\alpha-\{0\} \to \C^{n+1}-\{0\}$ is equivariant with respect to
$\C^*$, and since the quotient is $F : X \to \P^n$, one easily sees
that $F$ and $\widetilde F$ have the same degree.  This completes the
proof of the theorem.$\ \ \diamond$
\bigskip

\noi {\bf Remarks 5.9.}\ (i)\ Notice that in general, the integral
(5.7) is slightly different from the Grothendieck residue defined in
(4.8).  This is because $X_\alpha$ need not be simplicial at the point
$0 \in X_\alpha$.
\medskip
\noi (ii)\ When $X = \P^n$ and $F_0,\ldots,F_n$ are homogeneous of
degree $d$, note that the residue of Theorem 5.8 is computed {\it
not\/} on $\C^{n+1}$, but rather on $X_d = {\rm Spec}(\oplus_{k\ge0}
\C[x_0,\ldots,x_n]_{k\,d})$, which is the quotient of $\C^{n+1}$ by
the diagonal action of the $d$th roots of unity $\mu_d$.  Furthermore,
one can show that the Euler form of $X_d$ is $\Omega_d =
d\,dx_0\wedge\cdots\wedge dx_n$.

	Since $X_d$ is simplicial at the origin, the local residue
$\res_{0 \in {X_d}}(\omega_F(H))$ is defined, and combining Theorem 5.8
and equation (4.8), we see that
$$
\res_F(H) = \res_{0 \in X_d}\biggl({H\,\Omega_d\over F_0\cdots
F_n}\biggr) = \res_{0 \in \C^{n+1}}\biggl({H\,dx_0\wedge\cdots\wedge
dx_n\over F_0\cdots F_n}\biggr).
$$
Thus the toric residue equals both of the local residues that can be
defined in this situation, and Theorem 5.8 gives the toric
generalization of the first of these equalitites.
\bigskip
\bigskip

\noi {\bf References}
\parindent=.45truein
\bigskip
\item{[B]} W.~Baily, {\it The decomposition theorem for
$V$-manifolds\/}, Amer.~J.~Math.~{\bf 78} (1956), 862--888.
\medskip
\item{[BB]} V.~Batyrev and L.~Borisov, {\it Dual cones and mirror
symmetry for generalized Calabi-Yau manifolds\/}, to appear in {\sl
Mirror Symmetry II\/} (S.~T.~Yau, editor), alg-geom 9402002.
\medskip
\item{[BC]} V.~Batyrev and D.~Cox, {\it On the Hodge structure of
projective hypersurfaces in toric varieties\/}, Duke J.~Math. {\bf 75}
(1994), 293--338.
\medskip
\item{[BH]} W.~Bruns and J.~Herzog, {\sl Cohen-Macaulay Rings\/},
Cambridge Univ.~Press, Cambridge, 1993.
\medskip
\item{[CDS]}
E.~Cattani, A.~Dickenstein and B.~Sturmfels,
{\it Computing multidimensional residues\/}, to appear in
{\sl Proceedings MEGA94\/}, Birkh\"auser,  alg-geom 9404011.
\medskip
\item{[CH]}
N.~Coleff and M.~Herrera, {\sl Les Courants Residuels Associ\'es \`a
une Forme Meromorphe\/}, Lecture Notes in Math.~{\bf 633},
Springer-Verlag, Berlin Heidelberg New York, 1978.
\medskip
\item{[C1]} D.~Cox, {\it The homogeneous coordinate ring of a toric
variety\/}, J.~Algebraic Geom. {\bf 4} (1995), 17--50.
\medskip
\item{[C2]}
D.~Cox, {\it Toric  residues\/}, Ark.~Mat., to appear, alg-geom 9410017.
\medskip
\item{[Da]} V.~Danilov, {\it The geometry of toric varieties\/},
Russian Math.~Surveys {\bf 33} (1978), 97--154.
\medskip
\item {[Di]}
A.~Dickenstein, {\it Residues and ideals\/}, Aportaciones
Matem\'aticas, Notas de Investigaci\'on {\bf 5} (1992), 3--19.
\medskip
\item{[DK]} J.~Distler and S.~Kachru, {\it Duality of $(0,2)$ string
vacua\/}, hep-th 9501111.
\medskip
\item{[F]} W.~Fulton, {\sl Introduction to Toric Varieties\/},
Princeton Univ.~Press, Princeton, 1993.
\medskip
\item{[G]} P. Griffiths, {\it Variations on a theorem of Abel\/},
Inventiones Math. {\bf 35} (1976), 321--390.
\medskip
\item {[GH]}
P.~Griffiths and J.~Harris, {\sl Principles of Algebraic
Geometry\/}, John Wiley \& Sons, New York, 1978.
\medskip
\item {[H]}
R.~Hartshorne, {\sl Ample Subvarieties of Algebraic Varieties\/},
Lecture Notes in Math. {\bf 156}, Springer-Verlag, Berlin Heidelberg
New York, 1970.
\medskip
\item{[KSZ]} M.~Kapranov, B.~Sturmfels and A.~Zelevinsky, {\it
Quotients of toric varieties\/}, Math. Ann. {\bf 290} (1991),
643--655.
\medskip
\item{[K1]} A.~G.~Khovanskii, {\it Newton polyhedra and toroidal varieties\/},
Functional Anal. Appl. {\bf 11} (1977), 289--296.
\medskip
\item{[K2]} A.~G.~Khovanskii, {\it Newton polyhedra and the
Euler-Jacobi formula\/},
Russian Math. Surveys {\bf 33} (1978), 237--238.
\medskip
\item{[KK]} M.~Kreuzer and
E.~Kunz, {\it Traces in strict Frobenius algebras and strict complete
intersections\/}, J.~Reine angew.~Math. {\bf 381} (1987), 181--204.
\medskip
\item{[O]} T.~Oda, {\sl Convex Bodies and Algebraic Geometry\/},
Springer-Verlag, Berlin Heidelberg New York, 1988.
\medskip
\item {[PS]}
C.~Peters and J.~Steenbrink, {\it Infinitesimal variation of Hodge
structure and the generic Torelli theorem for projective
hypersurfaces\/}, in {\sl Classification of Algebraic and Analytic
Manifolds} (K.~Ueno, editor), Progress in Math.~{\bf 39},
Birkh\"auser, Boston Basel Berlin, 1983, 399--463.
\medskip
\item{[P]} D.~Prill, {\it Local classification of quotients
of complex manifolds by discontinuous groups\/},
Duke Math. J., {\bf 34} (1967), 375-386.
\medskip
\item{[R]} S.-S.~Roan, {\it Picard groups of hypersurfaces in toric
varieties\/}, preprint, 1995.
\medskip
\item{[Sa]} I.~Satake, {\it On a generalization of the notion of
manifolds\/}, Proc. Nat. Acad. Sci. U.S.A., {\bf 42} (1956), 359-363.
\medskip
\item{[St]} J.~Steenbrink, {\it Mixed Hodge structure on the vanishing
cohomology\/}, in {\sl Real and Complex Singularities: Proceedings Nordic
Summer School, Oslo 1976\/}, Sijthoff and Noordhoff, Alphen
aan den Rijn 1977, 525--563.
\medskip
\item {[Z]} G.~Ziegler, {\sl Lectures on Polytopes\/},
Springer-Verlag, Berlin Heidelberg New York, 1995.
\medskip
\bye